\documentclass[11pt,reqno]{amsart}
\usepackage{amsmath, amsfonts, amsthm, amssymb, psfig}
\usepackage{float,epsf,subfigure} 
\restylefloat{figure}

\theoremstyle{plain}
\newtheorem{theorem}{Theorem}[section]
\newtheorem{lemma}{Lemma}[section]
\newtheorem{proposition}{Proposition}[section]

\theoremstyle{definition}
\newtheorem{definition}{Definition}[section]

\theoremstyle{remark}
\newtheorem{remark}{Remark}[section]

\numberwithin{equation}{section}

\renewcommand{\r}{\rho}
\newcommand{\g}{\gamma}
\renewcommand{\t}{\theta}

\newcommand{\e}{\varepsilon}
\newcommand{\s}{\sigma}

\renewcommand{\div}{\text{div}}
\renewcommand{\O}{\Omega}
\renewcommand{\l}{\lambda}
\renewcommand{\d}{\partial}
\renewcommand{\a}{\alpha}

\newcommand{\R}{{\mathbb R}}

\def\be{\begin{equation}}
\def\ee{\end{equation}}
\def\bes{\begin{equation*}}
\def\ees{\end{equation*}}

\def\bc{\begin{cases}}
\def\ec{\end{cases}}

\begin{document}

\title[Vanishing Viscosity Method for Transonic Flow]
{Vanishing Viscosity Method for Transonic Flow}
\author[G.-Q. Chen]{Gui-Qiang Chen}
\address{G.-Q. Chen, Department of Mathematics, Northwestern University,
         Evanston, IL 60208.}
\email{\tt gqchen@math.northwestern.edu}

\author[M. Slemrod]{Marshall Slemrod}
\address{M. Slemrod, Department of Mathematics, University of Wisconsin,
Madison, WI 53706.} \email{\tt slemrod@math.wisc.edu}

\author[D. Wang]{Dehua Wang}
\address{D. Wang, Department of Mathematics, University of Pittsburgh,
                Pittsburgh, PA 15260.}
\email{\tt dwang@math.pitt.edu}

\dedicatory{Dedicated to Constantine M. Dafermos on the Occasion of
His 65th Birthday}

\date{\today}

\begin{abstract}
A vanishing viscosity method is formulated for two-dimensional
transonic steady irrotational compressible fluid flows with
adiabatic constant $\gamma\in [1,3)$. This formulation allows a
family of invariant regions in the phase plane for the corresponding
viscous problem, which implies an upper bound uniformly away from
cavitation for the viscous approximate velocity fields. Mathematical
entropy pairs are constructed through the Loewner-Morawetz relation
by entropy generators governed by a generalized Tricomi equation of
mixed elliptic-hyperbolic type, and the corresponding entropy
dissipation measures are analyzed so that the viscous approximate
solutions satisfy the compensated compactness framework. Then the
method of compensated compactness is applied to show that a sequence
of solutions to the artificial viscous problem, staying uniformly
away from stagnation, converges to an entropy solution of the
inviscid transonic flow problem.
\end{abstract}
\maketitle

\section{Introduction}

In two significant papers written a decade apart, Morawetz
\cite{Mor85,Mor95} presented a program for proving the existence of
weak solutions to the equations governing two-dimensional steady
irrotational inviscid  compressible flow in a channel or exterior to
an airfoil. As is well known, the classical results of Shiffman
\cite{Shif1} and Bers \cite{Bers-book} apply when the upstream speed
is sufficiently small, for which the flow remains subsonic and the
governing equations are elliptic (also see
\cite{Dong-book,Finn1,Finn2,FG1,FG3,Gilbarg,GSe,GSh}). However,
beyond a certain speed at infinity (determined by the flow
geometry), the flow becomes transonic which, coupled with
nonlinearity, yields shock formation (cf. \cite{Mor82}). Morawetz's
program in \cite{Mor85,Mor95} was to imbed the problem within an
assumed viscous framework for which the compensated compactness
framework (see Section 6 of this paper) would be satisfied. Under
this assumption, Morawetz
proved that solutions of the as yet unidentified viscous problem
have a convergent subsequence whose limit is a solution of the
transonic flow problem.

The purpose of this paper is to present such a viscous formulation,
hence completing part but not all of Morawetz's program.
Specifically, a vanishing viscosity method is formulated for
two-dimensional transonic steady irrotational compressible fluid
flows with adiabatic constant $\gamma\in [1, 3)$ to ensure a family
of invariant regions for the corresponding viscous problem, which
implies an upper bound uniformly away from cavitation for the
viscous approximate velocity fields. Mathematical entropy pairs are
constructed through the Loewner-Morawetz relation by entropy
generators governed by a generalized Tricomi equation of mixed
elliptic-hyperbolic type, and the corresponding entropy dissipation
measures are analyzed so that the viscous approximate solutions
satisfy the compensated compactness framework.
%
Then the method of compensated compactness is applied to show that a
sequence of solutions to the viscous problem, staying uniformly away
from stagnation, converges to an entropy solution of the inviscid
transonic flow problem.

On the other hand, Morawetz's assumption of no stagnation points for
the flow has not been removed; and the slip boundary condition
$(u,v)\cdot {\bf n}=0$ on the obstacle is only satisfied as the
inequality $(u,v)\cdot{\bf n}\ge 0$ for the general case (the
desired condition $(u,v)\cdot{\bf n}=0$ may be achieved when a
solution has no jump with certain regularity along the obstacle),
where $(u,v)$ is the fluid velocity field and ${\bf n}$ is the unit
normal on the obstacle pointing into the flow. Both issues may
reflect rather complicated boundary layer behavior for transonic
flow \cite{Rusak} and require further investigation, which are the
topics for subsequent research. In addition, we note that, when
$\gamma\ge 3$, cavitation is indeed possible and will be the topic
of a sequel to this paper.

This paper is divided into nine sections after the introduction. In
Section 2, the classical fluid equations are presented for
isentropic and isothermal irrotational planar flow. In Section 3,
the fluid equations are rewritten in the polar variables $(u,
v)=(q\cos\t, q\sin\t)$. In Section 4, we analyzes the behavior of
the Riemann invariants in the supersonic region, which guides the
design of artificial viscous terms in our vanishing viscosity
method. In Section 5, we continue our analysis and set boundary
conditions for the viscous system. This yields a family of invariant
regions in the $(u,v)$ fluid phase plane for $\gamma\in [1, 3)$,
where $\g$ is the ratio of specific heats if $\g>1$ and $\g=1$
denotes the case of constant temperature. In particular, any
invariant region is uniformly bounded away from the cavitation
circle in the phase plane, which yields a upper bound uniformly away
from cavitation in the viscous approximate solutions.
In Section 6, we formulate a compensated compactness framework for
steady flow. In Section 7, we first construct all mathematical
entropy pairs for the potential flow system
through the Loewner-Morawetz relation by entropy generators governed
by the generalized Tricomi equation of mixed elliptic-hyperbolic
type. Then we introduce the notion of entropy solutions through an
entropy pair by a convex entropy generator suggested by the work of
Osher-Hafez-Whitlow \cite{OHW}.
In Section 8, we use the previously derived uniform $L^\infty$
bounds on $(u^\epsilon,v^\epsilon)$, plus the convex entropy
generator introduced in Section 7 to guarantee our problem is within
the compensated compactness framework. In Section 9, we develop
Morawetz's argument in \cite{Mor95} for proving convergence of a
subsequence of solutions
of our viscous problem to the irrotational fluid problem.
Finally, in Section 10, we prove the existence of smooth solutions
to our viscous problem.

\section{Mathematical Equations}
Consider the artificial viscous system in a
domain $\O\subset\R^2$:
\begin{equation} \label{vcl}
\begin{cases}
v_x-u_y=R_1, \\
(\r u)_x+(\r v)_y=R_2,
\end{cases}
\end{equation}
where $R_1$ and $R_2$ are the artificial viscosity terms to be
determined, and $(u, v)$ is the flow velocity field. For a
polytropic gas with the adiabatic exponent $\g>1$,
$p=p(\r)=\r^\g/\g$ is the normalized pressure, the renormalized
density $\r$ is given by Bernoulli's law:
\begin{equation} \label{ber1}
\r=\r(q)=\big(1-\frac{\g-1}2q^2\big)^{\frac1{\g-1}},
\end{equation}
where $q$ is the flow speed defined by $q^2=u^2+v^2$. The sound
speed $c$ is defined as
\begin{equation} \label{sp}
c^2=p'(\r)=1-\frac{\g-1}2q^2.
\end{equation}
At the cavitation point $\r=0$,
 $$q=q_{cav}:= \sqrt{\frac2{\g-1}}.$$
 At the stagnation point $q=0$, the density reaches its maximum
 $\r=1$. Bernoulli's law \eqref{ber1} is valid for $0\le q\le
 q_{cav}$.  At the sonic point
$q=c$, \eqref{sp} implies $q^2=\frac2{\g+1}$. Define the critical
speed $q_{cr}$ as
$$q_{cr}:=\sqrt{\frac2{\g+1}}.$$
We rewrite Bernoulli's law \eqref{ber1} in the form
\begin{equation} \label{ber9}
q^2-q_{cr}^2=\frac2{\g+1}\left(q^2-c^2\right).
\end{equation}
 Thus the flow is subsonic when $q<q_{cr}$, sonic when $q=q_{cr}$,
 and supersonic when $q>q_{cr}$.

For the isothermal flow ($\g=1$), $p={c}^2\r$ where ${c}>0$ is the
constant sound speed, and the density $\r$ is given by Bernoulli's
law:
\begin{equation} \label{ber2}
\r=\r_0\exp\big(-\frac{u^2+v^2}{2{c}^2}\big)
\end{equation}
for some constant $\r_0>0$. In this case, $q_{cr}={c}$. After
scaling, we can take $c=\r_0=1$.

\section{Formulation in Polar Coordinate Phase Plane}

We now use the polar coordinates in the phase plane:
$$
u=q\cos \t, \qquad v=q\sin\t,
$$
and rewrite the viscous conservation laws \eqref{vcl} in terms of
$(q,\t)$. Write the second equation  of \eqref{vcl} as
$$
\r_x u+\r u_x+\r_y v+\r v_y=R_2
$$
or
$$
\r'(q)q_x u+\r u_x+\r'(q)q_y v+\r v_y=R_2,
$$
and use
$$
q=\sqrt{u^2+v^2}, \quad q_x=\frac1{q}(uu_x+vv_x), \quad
q_y=\frac1{q}(uu_y+vv_y), \quad \r'(q)=-\frac{\r q}{c^2}
$$
to find
$$
(c^2-u^2)u_x-uv(v_x+u_y)+(c^2-v^2)v_y=\frac{c^2}{\r}R_2.
$$
Then \eqref{vcl} becomes
$$
A\begin{bmatrix} u\\ v \end{bmatrix}_x
 +B\begin{bmatrix} u\\ v \end{bmatrix}_y=\begin{bmatrix} -R_1\\ \frac{c^2}{\r}R_2
 \end{bmatrix},
 $$
where
$$
A=\begin{bmatrix}
      0& -1\\
      c^2-u^2& -uv
      \end{bmatrix}, \qquad
   B=\begin{bmatrix}
     1 & 0\\
     -uv &c^2-v^2
     \end{bmatrix}.
$$
Thus, in terms of $(q,\t)$, we obtain
\begin{equation} \label{a1}
A_1\begin{bmatrix} q\\ \t
\end{bmatrix}_x
 +B_1\begin{bmatrix} q\\ \t \end{bmatrix}_y
 =\begin{bmatrix} -R_1\\ \frac{1}{\r q}R_2
 \end{bmatrix},
\end{equation}
where
$$
A_1=\begin{bmatrix}
      -\sin\t& -q\cos\t\\
      \frac{c^2-q^2}{c^2q}\cos\t& -\sin\t
      \end{bmatrix}, \qquad
   B_1=\begin{bmatrix}
     \cos\t & -q\sin\t \\
     \frac{c^2-q^2}{c^2q}\sin\t& \cos\t
     \end{bmatrix}.
$$
\bigskip

In terms of $(\r,\t)$, using
$$\r_x=\r_qq_x=-\frac{\r q}{c^2}q_x, \quad \r_y=-\frac{\r q}{c^2}
q_y,$$ and thus
$$\begin{bmatrix} q\\ \t \end{bmatrix}_x
=\begin{bmatrix} -\frac{c^2}{\r q}\r_x\\ \t_x \end{bmatrix}, \qquad
\begin{bmatrix} q\\ \t \end{bmatrix}_y
=\begin{bmatrix} -\frac{c^2}{\r q}\r_y\\ \t_y \end{bmatrix},$$
 we find
\begin{equation} \label{a2}
A_2\begin{bmatrix} \r \\ \t \end{bmatrix}_x
 +B_2\begin{bmatrix} \r \\ \t \end{bmatrix}_y
 =\begin{bmatrix} -R_1\\ -R_2 \end{bmatrix},
\end{equation}
where
$$
A_2=\begin{bmatrix}
      \frac{c^2}{\r q}\sin\t& -q\cos\t\\
      \frac{1}{q}(c^2-q^2)\cos\t& \r q\sin\t
      \end{bmatrix}, \qquad
   B_2=\begin{bmatrix}
     -\frac{c^2}{\r q}\cos\t & -q\sin\t \\
     \frac{1}{q}(c^2-q^2)\sin\t&  -\r q\cos\t
     \end{bmatrix}.
$$

We symmetrize \eqref{a2} to obtain
\begin{equation} \label{a3}
A_3\begin{bmatrix} \r \\ \t \end{bmatrix}_x
 +B_3\begin{bmatrix} \r \\ \t \end{bmatrix}_y
 =\begin{bmatrix} -R_1\\ \frac{q^2}{c^2-q^2}R_2 \end{bmatrix},
\end{equation}
 where
$$
 A_3=\begin{bmatrix}
      \frac{c^2}{\r q}\sin\t& -q\cos\t\\
      -q\cos\t &\frac{-\r q^3\sin\t}{c^2-q^2}
      \end{bmatrix}, \qquad
   B_3=\begin{bmatrix}
     -\frac{c^2}{\r q}\cos\t & -q\sin\t \\
     -q\sin\t &\frac{\r q^3\cos\t}{c^2-q^2}
       \end{bmatrix}.
$$

\section{Choice of Artificial Viscosity from an Analysis of
           Riemann Invariants}

The two matrices $A_1, B_1$ in \eqref{a1} commute, thus their
transposes commute and they have common eigenvectors. The
eigenvalues of $A_1$ and $B_1$ are
$$
\l_\pm=-\sin\t\pm \frac{\sqrt{q^2-c^2}}{c}\cos\t,\quad
\mu_\pm=\cos\t\pm \frac{\sqrt{q^2-c^2}}{c}\sin\t.
$$
The left eigenvectors of $A_1$ are
$$
(\mp\frac{\sqrt{q^2-c^2}}{qc},1),
$$
and thus the Riemann invariants $W_\pm$ satisfy
\be\label{R}
 \frac{\d W_\pm}{\d\t}=1,
\quad \frac{\d W_\pm}{\d q}
 =\mp\frac{\sqrt{q^2-c^2}}{qc}
 \qquad\text{for}\,\, q\ge c.
 \ee
Multiply \eqref{a1} by $(\frac{\d W_\pm}{\d q},\frac{\d W_\pm}{\d
\t})$ to obtain
\begin{equation} \label{W}
\l_\pm\frac{\d W_\pm}{\d x}+\mu_\pm\frac{\d W_\pm}{\d y} =-\frac{\d
W_\pm}{\d q}R_1+\frac{1}{\r q}\frac{\d W_\pm}{\d\t}R_2.
\end{equation}

We now consider the viscosity terms of the form:
\begin{equation} \label{R2}
R_1=\e\nabla\cdot\left(\s_1(\r,\t)\nabla\t\right), \quad
  R_2=\e\nabla\cdot\left(\s_2(\r,\t)\nabla\r\right).
  \end{equation}
In particular, for a special choice of $\s_1$ and $\s_2$, we have
the desired relations for the Riemann invariants.

\begin{proposition} \label{LW}
If
\begin{equation} \label{s20}
\s_1=1, \quad \s_2=1-\frac{c^2}{q^2} \qquad\text{for}\,\,\, q>c,
  \end{equation}
then the Riemann invariants $W_\pm$ satisfy the following equations:
\begin{equation} \label{R20}
 \frac{qc}{\sqrt{q^2-c^2}}
 \left(\l_\pm\frac{\d W_\pm}{\d x}+\mu_\pm\frac{\d W_\pm}{\d y}\right)
 \mp\e\Delta W_\pm
 =-\e c \frac{(\g-3)q^2+4c^2}{2\r^2q^2\sqrt{q^2-c^2}}
  |\nabla\r|^2
\end{equation}
for $\gamma\ge 1$.
 \end{proposition}

\begin{proof} We first focus on the ``+"  part, since the ``-'' part can
be done analogously.
The proof is divided into three steps.
\medskip

\noindent  {\em Step 1:}\quad Substitute \eqref{R} into \eqref{W} to
obtain
\begin{equation}\label{4.6a}
 \l_+\frac{\d W_+}{\d x}+\mu_+\frac{\d
W_+}{\d y} =\frac{\sqrt{q^2-c^2}}{qc}R_1 +\frac{1}{\r q}R_2.
\end{equation}
Multiplication of \eqref{4.6a} by $\frac{qc}{\sqrt{q^2-c^2}}$ gives
$$
\frac{qc}{\sqrt{q^2-c^2}}\left(\l_+\frac{\d W_+}{\d x}
 +\mu_+\frac{\d W_+}{\d y}\right)
  =\frac{\d W_+}{\d\t}R_1 +\frac{c}{\r \sqrt{q^2-c^2}}R_2.
$$
Using $\frac{d\r}{dq}=-\r q/c^2$, we have
$$
\frac{\d W_+}{\d\r}=\frac{\d W_+}{\d q}\frac{dq}{d\r}
=\frac{c\sqrt{q^2-c^2}}{\r q^2},
$$
 and then
$$
\frac{c}{\r \sqrt{q^2-c^2}}
  =\frac{q^2}{q^2-c^2}\frac{\d W_+}{\d\r}.
$$
Therefore,
\begin{equation} \label{W2}
\frac{qc}{\sqrt{q^2-c^2}}\left(\l_+\frac{\d W_+}{\d x}
 +\mu_+\frac{\d W_+}{\d y}\right)
  =\frac{\d W_+}{\d\t}R_1 +\frac{q^2}{q^2-c^2}\frac{\d W_+}{\d\r}R_2.
\end{equation}
For the choice of viscosity terms in \eqref{R2}, we obtain
\be
\label{W3}
\begin{split}
&\frac{qc}{\sqrt{q^2-c^2}}
\left(\l_+\frac{\d W_+}{\d x}+\mu_+\frac{\d W_+}{\d y}\right)\\
&=\e\frac{\d W_+}{\d\t}\nabla\cdot\left(\s_1(\r,\t)\nabla\t\right)
 +\e \frac{q^2}{q^2-c^2}\frac{\d W_+}{\d \r}\nabla\cdot\left(\s_2(\r,\t)\nabla\r\right)\\
&=\e\frac{\d
U_+}{\d\t}\left(\frac{\d\s_1}{\d\r}\nabla\r\cdot\nabla\t
      +\frac{\d\s_1}{\d\t}|\nabla\t|^2\right)
      +\e\s_1(\r,\t)\frac{\d W_+}{\d\t}\Delta\t \\
&\quad +\e \frac{q^2}{q^2-c^2}\frac{\d
W_+}{\d\r}\left(\frac{\d\s_2}{\d\r}|\nabla\r|^2
 +\frac{\d\s_2}{\d\t}\nabla\r\cdot\nabla\t\right)
 +\e\s_2(\r,\t)\frac{q^2}{q^2-c^2}\frac{\d W_+}{\d\r}\Delta\r.
\end{split}
\ee

\noindent  {\em Step 2:} \quad
From a direct calculation, we know
$$
\Delta W_\pm=\frac{\d^2 W_{\pm}}{\d\r^2}|\nabla\r|^2
+\frac{\d W_{\pm}}{\d\r}\Delta\r+\Delta\t,
$$
then substitution for $\Delta\t$ in \eqref{W3} gives \bes
\begin{split}
&\frac{qc}{\sqrt{q^2-c^2}}
\left(\l_+\frac{\d W_+}{\d x}+\mu_+\frac{\d W_+}{\d y}\right)\\
&= \e\frac{\d
W_+}{\d\t}\left(\frac{\d\s_1}{\d\r}\nabla\r\cdot\nabla\t
      +\frac{\d\s_1}{\d\t}|\nabla\t|^2\right)
 +\e\s_1(\r,\t)\frac{\d W_+}{\d\t}
 \left(\Delta W_+-\frac{\d^2W_{+}}{\d\r^2}|\nabla\r|^2
  -\frac{\d W_{+}}{\d\r}\Delta\r\right) \\
&\quad +\e\frac{q^2}{q^2-c^2}\frac{\d W_+}{\d\r}
 \left(\frac{\d\s_2}{\d\r}|\nabla\r|^2
 +\frac{\d\s_2}{\d\t}\nabla\r\cdot\nabla\t\right)
 +\e\s_2(\r,\t)\frac{q^2}{q^2-c^2}\frac{\d W_+}{\d\r}\Delta\r.
\end{split}
\ees

To eliminate $\Delta\r$,  we choose
 \bes 
\s_2(\r)=\s_1(\r)\frac{q^2-c^2}{q^2},
 \ees
where $\s_1$ and $\s_2$ are independent of $\t$. Then we have
\bes   
\begin{split}
&\frac{qc}{\sqrt{q^2-c^2}}
\left(\l_+\frac{\d W_+}{\d x}+\mu_+\frac{\d W_+}{\d y}\right)\\
&=\e\frac{\d\s_1}{\d\r}\frac{\d W_+}{\d\t}\nabla\r\cdot\nabla\t
 +\e\s_1(\r)\frac{\d W_+}{\d\t}
 \left(\Delta W_+-\frac{\d^2W_{+}}{\d\r^2}|\nabla\r|^2\right)
  +\e\frac{\d\s_2}{\d\r}
  \frac{q^2}{q^2-c^2}\frac{\d W_+}{\d\r}|\nabla\r|^2.
\end{split}
\ees
Finally, write $\nabla\t=\nabla W_+ -\frac{\d
W_{+}}{\d\r}\nabla\r$ so that
\be\label{W5}
\begin{split}
&\frac{qc}{\sqrt{q^2-c^2}}
  \left(\l_+\frac{\d W_+}{\d x}+\mu_+\frac{\d W_+}{\d y}\right)
-\e\frac{\d\s_1}{\d\r}\frac{\d W_+}{\d\t}\nabla\r\cdot\nabla W_+
-\e\s_1 \frac{\d W_+}{\d\t}\Delta W_+\\
&=-\e\frac{\d\s_1}{\d\r} \frac{\d W_{+}}{\d\r}\frac{\d
W_+}{\d\t}|\nabla\r|^2 -\e\s_1\frac{\d^2W_+}{\d\r^2}\frac{\d
W_+}{\d\t}|\nabla\r|^2 +\e\frac{\d\s_2}{\d\r}
 \frac{q^2}{q^2-c^2}\frac{\d W_+}{\d\r}|\nabla\r|^2.
\end{split}
\ee A convenient choice of $\s_1$ and $\s_2$ is as in \eqref{s20}:
$$
\s_1=1, \qquad \s_2=\frac{q^2-c^2}{q^2}.
$$
Thus, using $\frac{\d W_+}{\d\t}=1$, we obtain
 \be \label{w10}
\frac{qc}{\sqrt{q^2-c^2}}
  \left(\l_+\frac{\d W_+}{\d x}+\mu_+\frac{\d W_+}{\d y}\right)
   -\e\Delta W_+ \\
=-\e\left(\frac{\d^2W_+}{\d\r^2}
 -\frac{d\s_2}{d\r}\frac{q^2}{q^2-c^2}\frac{\d W_+}{\d\r}\right)
|\nabla\r|^2.
\ee

\noindent  {\em Step 3:} \quad We now compute the term
 $$\frac{\d^2W_+}{\d\r^2}
 -\frac{d\s_2}{d\r}\frac{q^2}{q^2-c^2}\frac{\d W_+}{\d\r}$$
in \eqref{w10}.

First we check the simple isothermal case ($\g=1$) for which
$$
\r=e^{-\frac{q^2}2}, \quad c=1; \qquad \frac{d\r}{dq}=-\r q.
$$
Then
\bes
 \begin{split}
 &\frac{\d W_+}{\d q}=-\frac{\sqrt{q^2-1}}{q},\\
 &\frac{\d W_+}{\d \r}=\frac{\d W_+}{\d q}\frac{dq}{d\r}
 =\frac{\sqrt{q^2-1}}{\r q^2},\\
 &\frac{\d^2W_+}{\d\r^2}=\frac{-q^4+2 q^2-2}{\rho^2q^4\sqrt{q^2-1}},
 \end{split}
 \ees
 and
 $$
 \s_2=\frac{q^2-1}{q^2}=1-\frac1{q^2}, \quad
  \frac{d\s_2}{d\r}=\frac{d\s_2}{d q}\frac{d q}{d\r}
   =-\frac2{\r q^4}.
   $$
Thus,
$$
 \frac{\d^2W_+}{\d\r^2}
 -\frac{d\s_2}{d\r}\frac{q^2}{q^2-c^2}\frac{\d W_+}{\d\r}
 =-\frac{q^2-2}{\r^2q^2\sqrt{q^2-1}},
$$
and \eqref{w10} becomes, in the case $\gamma=1$, $c=1$, \be
\label{w11}
 \frac{q}{\sqrt{q^2-1}}
  \left(\l_+\frac{\d W_+}{\d x}+\mu_+\frac{\d W_+}{\d y}\right)
-\e\Delta W_+
 =\e\frac{q^2-2}{\r^2q^2\sqrt{q^2-1}}|\nabla\r|^2.
\ee

We now consider the case $\gamma>1$ for which
$$
c^2=1-\frac{\g-1}2q^2, \qquad
\s_2=\frac{q^2-c^2}{q^2}=\frac{\g+1}2-\frac1{q^2},
$$
$$\frac{\d W_+}{\d q}=-\frac{\sqrt{q^2-c^2}}{qc},
\qquad \frac{d\r}{dq}=-\frac{\r q}{c^2}.
$$
Then
\bes
 \begin{split}
 \frac{\d W_+}{\d \r}&=\frac{\d W_+}{\d q}\frac{d q}{d\r}
 =-\frac{\sqrt{q^2-c^2}}{qc}\left(-\frac{c^2}{\r q}\right)
 =\frac{c\sqrt{q^2-c^2}}{\rho q^2},\\
 \frac{d\s_2}{d\r}&=\frac{d\s_2}{d q}\frac{d q}{d\r}
  =-\frac{2c^2}{\r q^4},\\
 \frac{\d^2W_+}{\d\r^2}
  &=\frac{\d}{\d q}\left(\frac{\d W_+}{\d \r}\right)\frac{d q}{d\r}
  =-\frac{c^2}{\r q}\frac{\d}{\d q}
  \left(\frac{\d W_+}{\d\r}\right)
  =-\frac{c^2}{\r q}\frac{\d}{\d q}\left(\frac{c\sqrt{q^2-c^2}}{\r q^2}
   \right)\\
  &=-\frac{c\left((3-\g)q^4+2(\g-3)q^2c^2+4c^4\right)}{2\r^2q^4\sqrt{q^2-c^2}}.
\end{split}
 \ees
  Thus,
 \bes
 \begin{split}
 &\frac{\d^2W_+}{\d\r^2}
 -\frac{d\s_2}{d\r}\frac{q^2}{q^2-c^2}\frac{\d W_+}{\d\r}\\
 &=-\frac{c}{2\r^2q^4\sqrt{q^2-c^2}}
   \left((3-\g)q^4+2(\g-3)q^2c^2+4c^4\right)+\frac{2c^2}{\r q^4}\frac{q^2}{q^2-c^2}
    \frac{c\sqrt{q^2-c^2}}{\r q^2}\\
 &=\frac{c\left((\g-3)q^2+4c^2\right)}{2\r^2q^2\sqrt{q^2-c^2}},
\end{split}
 \ees
 and \eqref{w10} becomes, in the case $\gamma>1$,
\be \label{w12}
\begin{split}
 \frac{qc}{\sqrt{q^2-c^2}}
  \left(\l_+\frac{\d W_+}{\d x}+\mu_+\frac{\d W_+}{\d y}\right)
-\e\Delta W_+
 =-\e c \frac{(\g-3)q^2+4c^2}{2\r^2q^2\sqrt{q^2-c^2}}|\nabla\r|^2.
\end{split}
\ee This equality is consistent with the case $\g=1$, and this lemma
is proved for $W_+$. The computation for  $W_-$ is done analogously.
\end{proof}

\section{Boundary Conditions, Invariant Regions, and $L^\infty$ Bounds}

In this section, we consider the level sets of $W_\pm$ to assign
appropriate boundary conditions and identify a family of invariant
regions for the viscous problem \eqref{vcl} with
\eqref{R2}--\eqref{s20}.

First we discuss the isothermal case: $\gamma=1$ and $c=1$. From
Proposition \ref{LW}, we find that, when $q>\sqrt{2}$,
$W_+(x,y):=W_+(q(x,y),\theta(x,y))$ cannot have an interior minimum
since $\Delta W_+<0$. Similarly, $W_-(x,y):=W_-(q(x,y),\theta(x,y))$
cannot have an interior maximum since $\Delta W_->0$. If we assume
for the moment that $W_+$ cannot have a minimum on the boundary of
our domain, then a solution of the viscous problem which crosses
from $q\le \sqrt{2}$ to $q>\sqrt{2}$ at $(x_0,y_0)$ must satisfy
$$W_+(x,y)\ge W_+(x_0,y_0)|_{q=\sqrt{2}}.$$
Similarly, if we assume for the moment that $W_-$ cannot have a
maximum on the boundary of our domain, then a solution of the
viscous problem which crosses from $q\le \sqrt{2}$ to $q>\sqrt{2}$
at $(x_0,y_0)$ must satisfy
$$
W_-(x,y)\le W_-(x_0,y_0)|_{q=\sqrt{2}}.
$$
From the definition of $W_\pm$, we have
$$
\frac{\d W_\pm}{\d q}=\mp\frac{\sqrt{q^2-1}}{q}.
$$
Using the substitution $q=\sec t$, we can integrate to find
$$
W_\pm=\t\mp\sqrt{q^2-1}\pm\arccos(q^{-1})+C_\pm
=\t\mp\left(\sqrt{q^2-1}-1\right)\pm\left(\arccos(q^{-1})-\frac\pi{4}\right),
$$
where the constants $C_\pm$ is set to be
$C_\pm=(\pm\sqrt{q^2-1}\mp\arccos(q^{-1}))|_{q=\sqrt{2}}
 =\pm(1-\frac\pi{4})$.

Along the level curves $W_\pm=const.$, we have
$$
\frac{d\t}{d q}=\pm\frac{\sqrt{q^2-1}}{q}.
$$
Thus, on the level set of $W_+$, $\frac{d\t}{d q}>0$ and $\t$ is
increasing when $q$ is increasing; while, on the level set of $W_-$,
$\frac{d\t}{d q}<0$ and $\t$ is decreasing when $q$ is
 increasing.

\begin{figure}[ht]
\centerline{\psfig{file=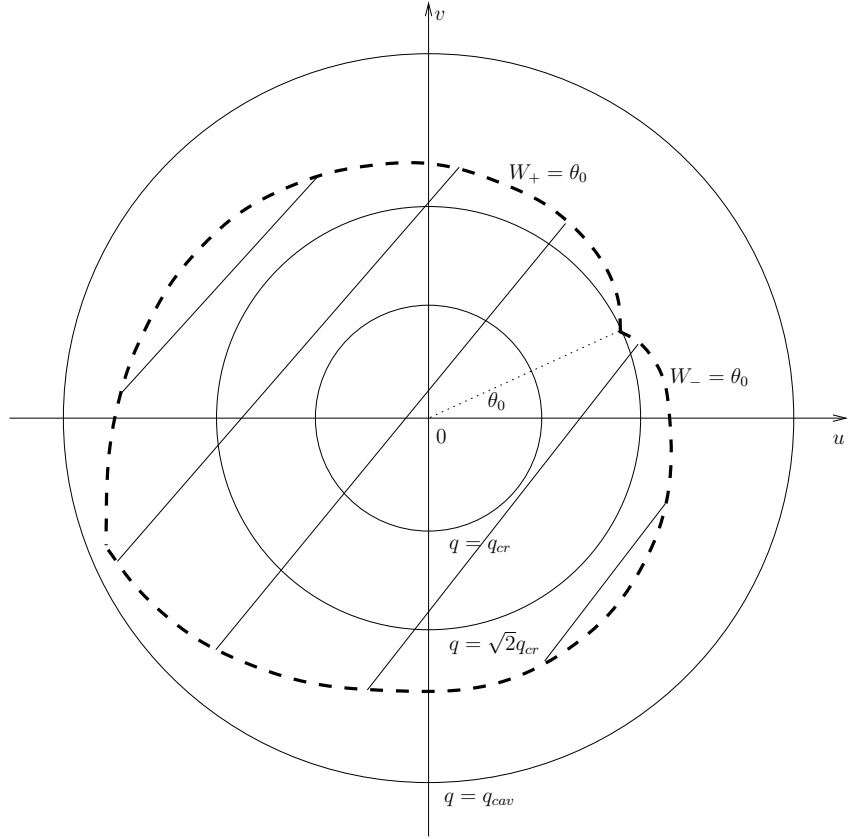,width=4.5in}} 
\caption{Invariant regions of ``apple" shape.}
 \label{ri}
\end{figure}

If $\t(x_0,y_0)=\t_0$, $q(x_0,y_0)=\sqrt{2}$ and we
 leave $q=\sqrt{2}$, we have
$$
\t(x,y)-\t_0\ge\left(\sqrt{q^2-1}-1\right)-\left(\arccos(q^{-1})-\frac\pi{4}\right)
 \qquad\,\,\text{(stay below $W_+$)},
$$
$$
\t(x,y)-\t_0\le
-\left(\sqrt{q^2-1}-1\right)+\left(\arccos(q^{-1})-\frac\pi{4}\right)
 \qquad\,\,\text{(stay above $W_-$)},
$$
that is,
$$
\left(\sqrt{q^2-1}-1\right)-\left(\arccos(q^{-1})-\frac\pi{4}\right)\le
|\t(x,y)-\t_0| \qquad\text{inside the ``apple" shaped region}.
$$
See Figure \ref{ri} with $q_{cav}=\infty$ when $\gamma=1$.

The same situation occurs when the level set curves
$W_\pm=W_\pm(q_0,\theta_0)$ for $\t(x_0,y_0)=\t_0$ and
$q(x_0,y_0)=q_0>\sqrt{2}$, for which we obtain similar invariant
regions of ``apple'' shape past the point $(q_0,\theta_0)$ in the
$(u,v)$-plane.

\medskip
Now we return to the issue of boundary conditions. Denote the
boundary of the bounded domain $\O$  by $\d\O$, the boundary of the
obstacle by $\d\O_1$, and the far field boundary by $\d\O_2$. Thus,
$\d\O=\d\O_1\cup\d\O_2$. Since we do not want $W_+$ to have a
minimum on boundary $\d\O_1$ and $W_-$ to have a maximum on boundary
$\d\O_1$, we  require
$$
\frac{\d W_+}{\d{\bf n}}<0, \quad
  \frac{\d W_-}{\d{\bf n}}>0
  \qquad\text{at all boundary points},
$$
where ${\bf n}$ is the unit normal into the flow region on
$\partial\Omega$. Recall
$$
\nabla W_\pm=\nabla\t\mp\frac{q^2-1}{q}\nabla q.
$$
Therefore, if we set
$$
\nabla\t\cdot{\bf n}=0 \qquad\text{on}\quad \d\O_1,
$$
then
$$
\text{sign}(\nabla W_\pm\cdot{\bf n}) =\mp\; \text{sign}(\nabla
q\cdot {\bf n})
$$
at those boundary points where $q>1$, i.e. where the Riemann
invariants are defined.
Hence, one resolution of the boundary
  condition issue is to set
$$
\e\s_2\nabla\r\cdot{\bf n}=-|\r\;(u,v)\cdot{\bf n}|
  \qquad\text{on}\quad \d\O_1,
$$
  and
$$
(u,v)-(u_\infty,v_\infty)=0  \qquad\text{on}\quad \d\O_2,
$$
with $q_\infty=|(u_\infty,v_\infty)|\in (0, q_{cav})=(0,\infty)$.
Trivially, $W_\pm$ are not even defined on $\d\O_2$ and hence cannot
  have a maximum or minimum there.
Since
$$
\nabla\r=\nabla q\frac{d\r}{d q}=-\r q\nabla q,
$$
we have automatically
$$
\text{sign}(\nabla W_\pm\cdot{\bf n}) =\pm \text{sign}(\nabla\r\cdot
{\bf n})=\mp 1 \qquad\mbox{on all boundaries},
$$
and our minimum principle for $W_+$ and maximum principle for $W_-$
are indeed valid.

Secondly, they formally yield
$$
\r\;(u,v)\cdot{\bf n}=0 \qquad\text{on}\quad \d\O_1,
$$
and
$$
(u,v)-(u_\infty,v_\infty)\to 0 \qquad\text{as}\quad
x^2+y^2\to\infty,
$$
if  $\e\to 0$. This is the case when any shock strength is zero at
its intersection point with the boundary as conjectured for the
shock formed in the supersonic bubble near the obstacle (see
Morawetz \cite{Mor04}). Finally, we see that the relation between
$\s_1$ and $\s_2$ is only necessary for $q\ge\sqrt{2}$. Any smooth
continuations of $\s_1$ and $\s_2$ inside $q=\sqrt{2}$ with
$\s_1,\s_2>0$ suffices.

With this, we conclude that $(q,\t)$ stay inside the ``apple" region
for $\gamma=1$. More generally, for $1<\gamma<3$, we have

\begin{theorem} \label{uniformbound}
Consider the viscous problem
\begin{equation} \label{vc}
\begin{cases}
v_x-u_y=\e\nabla\cdot\left(\s_1(\r)\nabla\t\right), \\
(\r u)_x+(\r v)_y=\e\nabla\cdot\left(\s_2(\r)\nabla\r\right),
\end{cases}
\end{equation}
with
$$
\s_1=1, \qquad \s_2=\frac{q^2-c^2}{q^2}
 \quad\text{for}\,\,\, q\ge \sqrt{2}q_{cr},
 \,\, 1\le\g <3,
$$
and the boundary conditions:
\begin{equation} \label{bc}
\begin{cases}
\nabla\t\cdot{\bf n}=0 \qquad\text{on}\quad \d\O_1,\\
\e\s_2\nabla\r\cdot{\bf n}=-|\r\;(u,v)\cdot{\bf n}|
  \qquad\text{on}\quad \d\O_1, \\
(u,v)-(u_\infty,v_\infty)=0  \qquad\text{on}\quad \d\O_2
 \,\,\text{with}\quad q_\infty<q_{cav},
\end{cases}
\end{equation}
Then all solutions to \eqref{vc}--\eqref{bc} are bounded on the
domain $\O$ staying uniformly in $\e$ away from cavitation, i.e.
there exists $q^*<q_{cav}$ such that $q^\varepsilon\le q^*$ that is
equivalent to $\rho(q^\varepsilon)\ge\underline{\rho}>0$ for some
$\underline{\rho}=\underline{\rho}(q^*)>0$. More specifically, when
$1\le \gamma<3$,
the viscous approximate solutions stay in a family of ``apple''
shaped invariant regions as shown in Figures {\rm 1} and {\rm 3--5}.
\end{theorem}

\begin{proof}
The theorem has been proved in the above for the isothermal case
$\g=1, c=1$, and $q_{cr}=1$.

For the case $1<\g<3$, from \eqref{R20}, we first require
$$
(\g-3)q^2+4c^2<0,
$$
which implies $q^2>\frac{4}{3-\g}c^2.$ Then, from Bernoulli's law
\eqref{ber9}, we have
$$
q^2-q_{cr}^2>\frac2{3-\g}c^2
 =\frac2{3-\g}\big(1-\frac{\g-1}2q^2\big),
$$
thus $q^2>2q_{cr}^2$, i.e. $q>\sqrt{2}q_{cr}$.

\begin{figure}[ht]
\centerline{\psfig{file=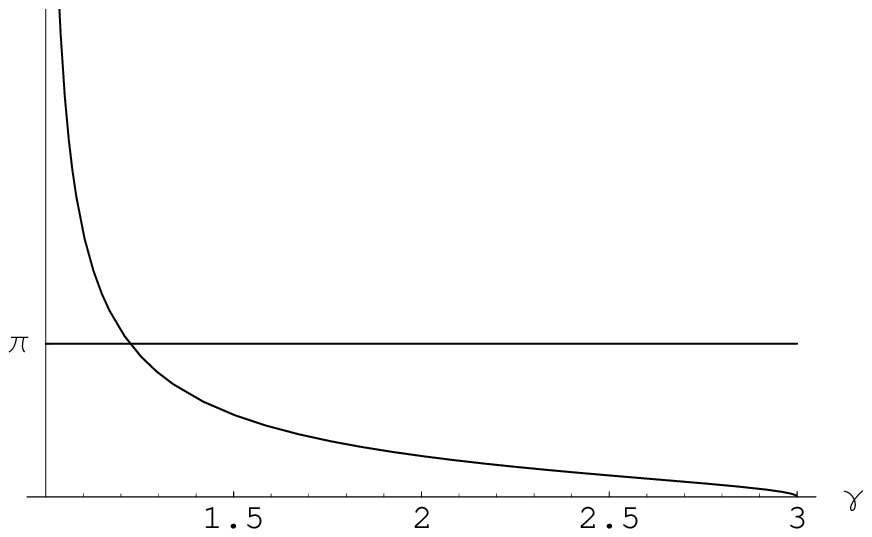,width=4in,height=3in}} 
\caption{The graph of $a(\g)$}
 \label{ag}
\end{figure}

\medskip
We first consider the case that the far field speed $q_\infty\le
\sqrt{2}q_{cr}$. Then, from the definition of the Riemann invariants
\eqref{R}, we follow Landau-Lifshitz \cite{LL-FM}, page 446, and
have
$$
W_\pm=\t\mp \left(W(q)-W(\sqrt2 q_{cr})\right),
$$
where
$$
W(q)=\sqrt\frac{\gamma+1}{\gamma-1}
\arcsin\sqrt{\frac{\gamma-1}2\Big(\frac{q^2}{q_{cr}^2}-1\Big)}
-\arcsin\sqrt{\frac{\gamma+1}2\Big(1-\frac{q_{cr}^2}{q^2}\Big)}.
$$

\begin{figure}[ht]
\centerline{\psfig{file=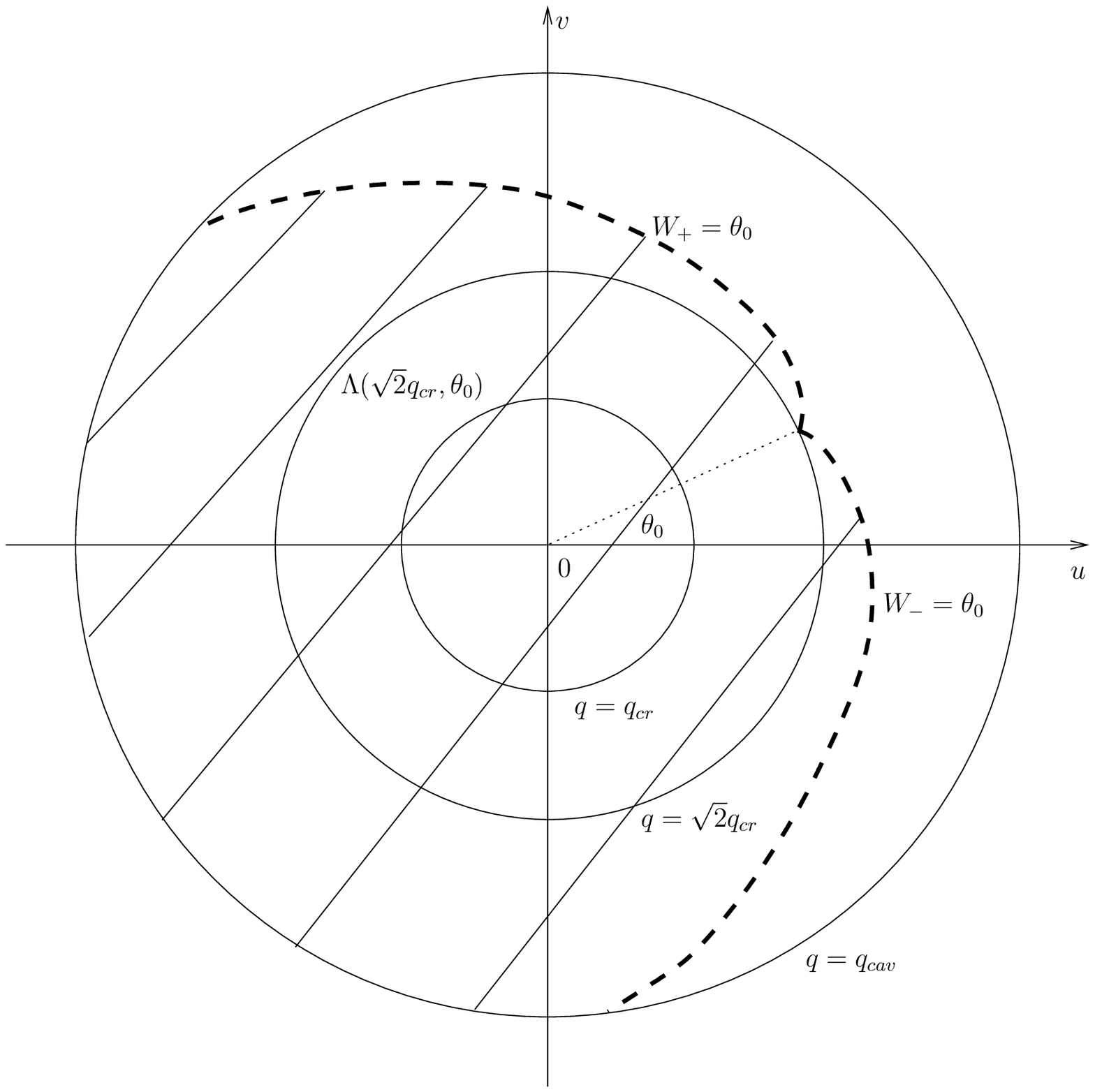,width=4in}} 
\caption{Round ``apple" shaped region $\Lambda(\sqrt{2}q_{cr},
\theta_0)$ when $a(\gamma)<\pi$}
 \label{ras}
\end{figure}

The two level set curves of $W_\pm$ intersect to form a invariant
region of ``apple" shape past $(\sqrt{2}q_{cr}, \theta_0)$ as in
Figure \ref{ri} if
$$
W(q_{cav})-W(\sqrt2 q_{cr})> \pi,
$$
that is,
\begin{equation*}
\begin{split}
 a(\g)&:=W(q_{cav})-W(\sqrt2 q_{cr}) \\
    & =\Big(\sqrt\frac{\gamma+1}{\gamma-1} -1\Big)\frac\pi 2
-\Big(\sqrt\frac{\gamma+1}{\gamma-1} \arcsin\sqrt{\frac{\gamma-1}2}
-\arcsin\sqrt{\frac{\gamma+1}4}\Big)> \pi.
\end{split}
\end{equation*}
 %

The graph of the function $a(\g)$ is shown in Figure \ref{ag}, which
shows that $a(\g*)=\pi$ for some $\g*\approx 1.224$ and $a(\g)>\pi$
for $1\le\g<\g*$.

\begin{figure}[ht]
\centerline{\psfig{file=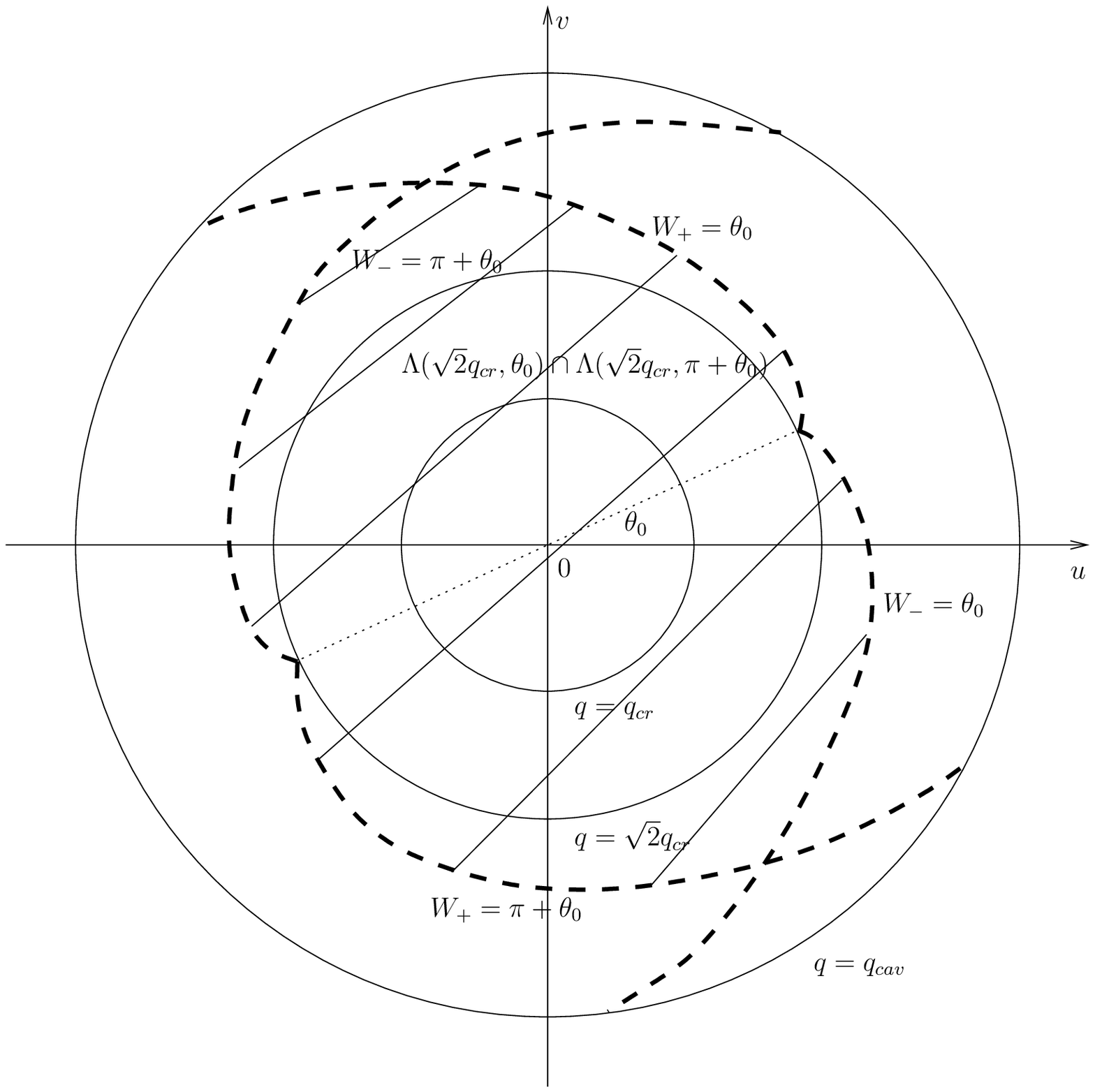,width=4in}} 
\caption{Invariant regions $\Lambda(\sqrt{2}q_{cr},\theta_0)\cap
\Lambda(\sqrt{2}q_{cr}, \pi+\theta_0)$ when
$a(\gamma)\in(\frac{\pi}{2}, \pi)$}
 \label{ras}
\end{figure}

For general $\g\in [\g^*, 3)$, the level set curves of
$W_\pm=\theta_0$ past $(\sqrt{2}q_{cr}, \theta_0)$ end at some
points on the cavitation circle $q=q_{cav}$ and do not intersect
each other. We denote $\Lambda(\sqrt{2}q_{cr},\theta_0)$ the round
``apple'' shaped region formed by $W_+>\theta_0$, $W_-<\theta_0$,
and the cavitation circle $q=q_{cav}$; see Figure 3.

\begin{figure}[ht]
\centerline{\psfig{file=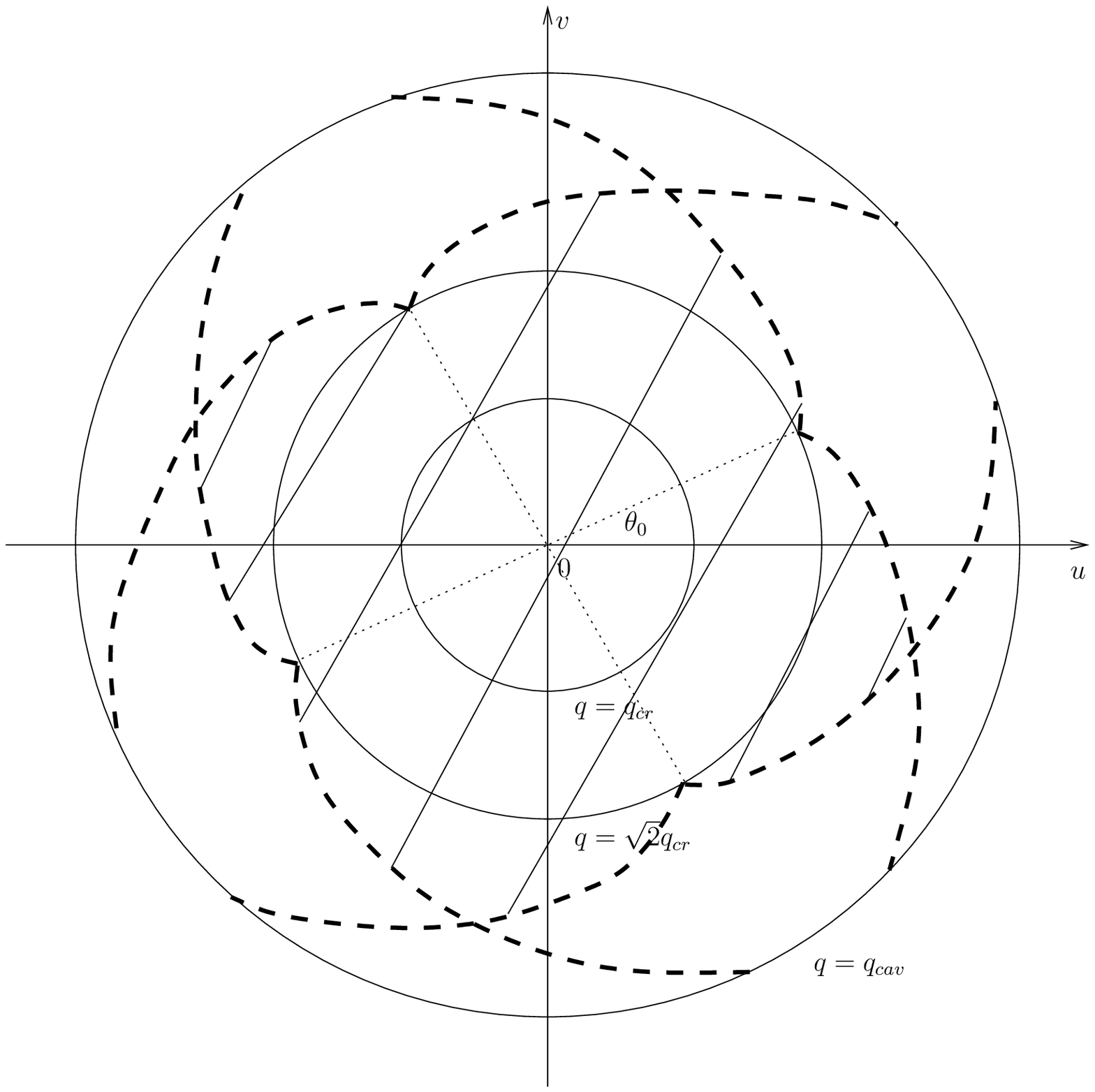,width=4in}} 
\caption{Invariant regions  $\Lambda(\sqrt{2}q_{cr},\theta_0)\cap
\Lambda(\sqrt{2}q_{cr},\frac{\pi}{2}+\theta_0)
\cap\Lambda(\sqrt{2}q_{cr},\pi+\theta_0)\cap\Lambda(\sqrt{2}q_{cr},\frac{3\pi}{2}+\theta_0)$
when $a(\gamma)\in (\frac{\pi}{4}, \frac{\pi}{2}]$}
 \label{ras}
\end{figure}

When $a(\gamma)\in (\frac{\pi}{2}, \pi]$, we find that, for any
$\theta_0 \in [0, 2\pi)$,
$$
\Lambda(\sqrt{2}q_{cr}, \theta_0)\cap \Lambda(\sqrt{2}q_{cr},
\pi+\theta_0)
$$
forms an invariant region for the viscous solutions; this follows
from application of our minimum and maximum principle for $W_+$ and
$W-$, respectively, at the crossing point on the new level set
curves. Such invariant regions stay away from the cavitation circle
(see Figure 4).

When $a(\gamma)\in (\frac{\pi}{4}, \frac{\pi}{2}]$, then, for any
$\theta_0 \in [0, 2\pi)$,
$$
\Lambda(\sqrt{2}q_{cr},\theta_0)\cap \Lambda(\sqrt{2}q_{cr},
\frac{\pi}{2}+\theta_0) \cap\Lambda(\sqrt{2}q_{cr},\pi+\theta_0)\cap
\Lambda(\sqrt{2}q_{cr},\frac{3\pi}{2}+\theta_0),
$$
forms an invariant region for the viscous solutions, staying away
from cavitation; see Figure 5.

In general, when $a(\gamma)\in (\frac{\pi}{2^n},
\frac{\pi}{2^{n-1}})$, $n=1,2,\dots$, for each $\theta_0\in
[0,2\pi)$, it requires an intersection of at least $2n$ round
``apple'' shaped regions including $\Lambda(\sqrt{2}q_{cr},
\theta_0)$ to form an invariant region staying away from the
cavitation circle.

\medskip
When the far field speed $q_\infty\in (\sqrt{2}q_{cr}, q_{cav})$,
then the level set curves of $W_\pm=W_\pm(q_0, \theta_0)$, for some
$q_0\in (q_\infty, q_{cav})$ and
$\theta_0=arctan\big(\frac{v_\infty}{u_\infty}\big)$, past the point
$(q_0, \theta_0)$ either intersects each other or end at some points
on the cavitation circle $q=q_{cav}$. As before, we denote
$\Lambda(q_0,\theta_0)$ the ``apple'' shaped region formed by
$W_+>W_+(q_0,\theta_0), W_-<W_-(q_0,\theta_0)$, and the cavitation
circle $q=q_{cav}$ similar to either Figure 1 or Figure 3. Then we
can similarly obtain the invariant regions which consist of either a
single ``apple'' shaped region or some unions of certain number
round ``apple'' shaped regions including $\Lambda(q_0,\theta_0)$,
which stay away from the cavitation but include the state
$(u_\infty, v_\infty)$ (cf. Figures 5--6).
\end{proof}

\section{Compensated Compactness Framework for Steady Flow}

Let a sequence of functions $w^\varepsilon(x,y)=(u^\varepsilon,
v^\varepsilon)(x,y)$, defined on open subset $\Omega\subset \R^2$,
satisfy the following Set of Conditions (A):

\medskip
{\rm (A.1)} $q^\varepsilon(x,y)=|w^\varepsilon(x,y)|\le q_{*}$\,\,
a.e. in $\Omega$, for some positive constant $q_*<q_{cav}<\infty$;

\smallskip
{\rm (A.2)} $\partial_xQ_{1\pm}(w^\varepsilon)+\partial_y
Q_{2\pm}(w^\varepsilon)$ are confined in a compact set in
$H_{loc}^{-1}(\O)$, for any entropy-entropy flux pairs $(Q_1,Q_2)$
so that $(Q_{1\pm}(w^\varepsilon), Q_{2\pm}(w^\varepsilon))$ are
confined in a bounded set uniformly in $L^\infty_{loc}(\Omega)$.
\medskip

\noindent Then, by the Div-Curl Lemma of Tartar \cite{Ta1} and Murat
\cite{Mu2} and the Young measure representation theorem for a
uniformly bounded sequence of functions (cf. Tartar \cite{Ta1}; also
Ball \cite{Ball}), we have the following commutation identity:
\begin{equation} \label{t1}
\begin{split}
&<\nu(w), Q_{1+}(w) Q_{2-}(w) -Q_{1-}(w) Q_{2+}(w)> \\
&\quad =<\nu(w), Q_{1+}(w)><\nu(w), Q_{2-}(w)>
        -<\nu(w), Q_{1-}(w)><\nu(w), Q_{2+}(w)>,
\end{split}
\end{equation}
where $\nu=\nu_{x,y}(w), w=(u,v),$ is the associated family of Young
measures (probability measures) for the sequence
$w^\varepsilon(x,y)=(u^\varepsilon, v^\varepsilon)(x,y)$. This is
equivalent to
\begin{equation} \label{t2}
<\nu(w)\otimes\nu(w'),\,\, I(w,w')>=0,
\end{equation}
where $\nu(w)\otimes\nu(w')$ is a product measure for
$(w,w')\in\R^2\times\R^2$ and
$$
I(w,w')=
 (Q_{1+}(w)-Q_{1+}(w'))(Q_{2-}(w)-Q_{2-}(w'))-(Q_{2+}(w)-Q_{2+}(w'))
 (Q_{1-}(w)-Q_{1-}(w')).
$$

The main point  for the compensated compactness framework is to
prove that $\nu$ is in fact a Dirac measure by using entropy pairs,
which implies the compactness of the sequence
$w^\varepsilon(x,y)=(u^\varepsilon, v^\varepsilon)(x,y)$ in
$L^1_{loc}(\Omega)$. Some additional references on compensated
compactness method include Chen \cite{Chen2,Chen3}, Dafermos
\cite{Dafermos-book}, DiPerna \cite{DiPerna2}, Evans \cite{evans},
and Serre \cite{Serre}. Theorem \ref{uniformbound} shows that (A.1)
is satisfied for our viscous solution sequence $w^\e=(u^\varepsilon,
v^\varepsilon)$. In the next two sections, we show how (A.2) can be
achieved.

\section{Entropy Pairs via Entropy Generators}

We  first construct all mathematical entropy pairs for the potential
flow system.

For some function $V(\r,\t)$ to be determined, multiply \eqref{a3}
from left by $(V_\t,V_\r)$ to get the  entropy equality
\be
\label{w30}
 Q_{1x}+Q_{2y}=-V_\t R_1+\frac{q^2}{c^2-q^2}V_\r R_2,
\ee
where $(Q_1,Q_2)$ is defined by \be \label{w31}
 \begin{split}
  &\frac{\d Q_1}{\d\r}
   =\frac{c^2}{\r q}\sin\t\; V_\t -q\cos\t\; V_\r, \quad
   \frac{\d Q_1}{\d\t}
   =-q\cos\t\; V_\t-\frac{\r q^3\sin\t}{c^2-q^2} V_\r, \\
  &\frac{\d Q_2}{\d\r}
   =-\frac{c^2}{\r q}\cos\t\; V_\t -q\sin\t\; V_\r, \quad
   \frac{\d Q_2}{\d\t}
   =-q\sin\t\; V_\t+\frac{\r q^3\cos\t}{c^2-q^2} V_\r.
   \end{split}
   \ee
We note that, the appearance of the term $c^2-q^2$ in the
   denominator of \eqref{w31} is only a consequence of our
   formalism. It will cancel out as we proceed.
Using $\frac{\d^2Q_i}{\d\t\d\r}=\frac{\d^2Q_i}{\d\r\d\t}$,
 $i=1,2$, from \eqref{w31}, we see that $V$ satisfies the Tricomi
type  equation of mixed type:
 \be \label{w32}
\frac{c^2}{\r q}V_{\t\t}+qV_\r
 +\left(\frac{\r q^3}{c^2-q^2}V_\r\right)_\r=0.
 \ee
Thus, we have defined an entropy pair $(Q_1,Q_2)$ generated by $V$.
Alternatively, we can define an entropy pair $(Q_1,Q_2)$ generated
by $H$, where $H$ and $V$ are related by
\begin{equation} \label{HV}
 \r H_{\mu\t}-H_\t=-V_\t, \qquad
   H_\mu+\frac1\r H_{\t\t}=\frac{q^2}{c^2-q^2}V_\r,
\end{equation}
with $\mu=\mu(\r)$  defined by $\mu'(\r)=c^2/q^2$, and
 $H$ is determined by the generalized Tricomi equation:
\begin{equation} \label{Tricomi}
H_{\mu\mu}+ \frac1{\r^2}(1-M^2)H_{\t\t}=0,
\end{equation}
and $M=q/c$ is the Mach number.

\begin{lemma} \label{Q8}
The entropy pairs $(Q_1,Q_2)$ are given by the Loewner-Morawetz
relation:
 \begin{equation} \label{Q}
 Q_1=\r q H_{\mu}\cos\t -q H_{\t}\sin\t, \qquad
 Q_2=\r q H_\mu\sin\t  +q H_\t\cos\t,
 \end{equation}
where the generators $H$ are all solutions of \eqref{Tricomi}.
\end{lemma}

This can be seen by differentiation of \eqref{Q} with respect to
$(\r,\t)$ and comparison with \eqref{w31}.

\medskip
A prototype of the generators $H$, as suggested in \cite{OHW}, is
\begin{equation}\label{h-entropy}
H^*=\frac{\t^2}2+\int^\mu\int^{\mu'}\frac1{\r^2}(M^2-1)d\mu'd\mu,
\end{equation}
which is a trivial solution of the generalized Tricomi equation
\eqref{Tricomi}. Notice that $H^*$ is strictly convex in
$(\mu,\theta)$ in the supersonic region. We denote by $(Q_1^*,
Q_2^*)$ the corresponding  entropy pair generated by the convex
generator $H^*$. With this, we introduce the notion of entropy
solutions.

\begin{definition}[Notion of Entropy Solutions]\label{def-entropysolution}
A bounded, measurable vector function $w(x,y):=(u,v)(x,y), q(x,y)\le
q_{cav}$, is called an entropy solution of the potential flow system
in a domain $\Omega$:
\begin{equation}\label{pfs}
\begin{cases}
v_x-u_y=0, \\
(\r u)_x+(\r v)_y=0,
\end{cases}
\end{equation}
if $w(x,y)$ satisfies \eqref{pfs} and the following entropy
inequality:
\begin{equation}\label{entropyinq}
Q^*_{1x}+Q^*_{2y}\le 0
\end{equation}
in the sense of distributions in $\Omega$, in addition to the
corresponding boundary conditions in the trace or asymptotic sense
on $\partial\Omega$.
\end{definition}

The physical correctness of the entropy inequality
\eqref{entropyinq} is provided by Theorem 2.1 of Osher-Hafez-Whitlow
\cite{OHW}.

\section{$H^{-1}$ Compactness of Entropy Dissipation Measures}

\begin{figure}[ht]
\centerline{\psfig{file=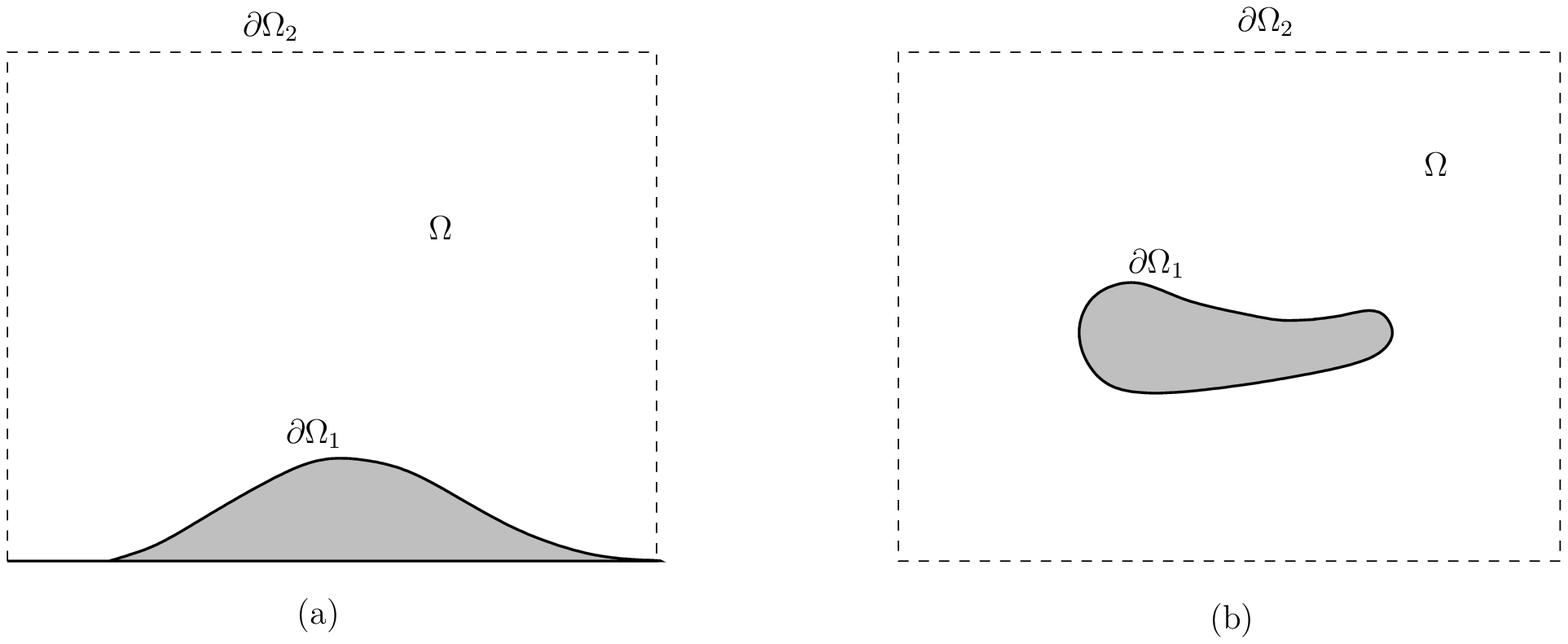,width=5.5in}} 
\caption{Domains}
 \label{domain}
\end{figure}

We now limit ourselves to the case $v_\infty=0$ and the two types of
domains $\O$ in Figure \ref{domain}(a)-(b), where $\d\O_1$ (the
solid curve in (a) and the solid closed curve in (b)) is the
boundary of the obstacle, $\d\O_2$ (dashed line segments in both (a)
and (b)) is the far field boundary, and $\O$ is the domain bounded
by $\d\O_1$ and $\d\O_2$. This assumption implies $\t=0$ on
$\partial\O_2$. In case (b), the circulation about the boundary
$\d\O_2$ is zero.

\begin{proposition} \label{I2}
Let $w^\e=(u^\e,v^\e)$ be a solution to \eqref{vc}--\eqref{bc} with
$u_\infty>0$ and $v_\infty=0$ on either domain $\Omega$ in Figure
{\rm 6(a)-(b)}, satisfying that $q^\e\le q^*<q_{cav}$ and
$\theta^\e$ is bounded. Then the integral
$$
\e\int_\O\left(\s_1(\r^\e)|\nabla\t^\e|^2
 +\s_2(\r^\e)\frac{c^2(\r^\e)}{(\r^\e q^\e)^2}|\nabla\r^\e|^2\right)dxdy
$$
is bounded uniformly in all $\e>0$.
\end{proposition}

\begin{proof}
We choose the special generator $H^*$ in \eqref{h-entropy}
which yields $V^*$ of the form: \be \label{V4}
 V^*=\frac{\t^2}2+P(\r), \qquad
 V^*_\r=P'(\r)=\frac{c^2-q^2}{q^2}\int_{\bar{q}}^q\frac{dq}{\r q},
  \qquad V^*_\t=\t.
\ee
and \eqref{w30} becomes \be \label{w34}
 Q^*_{1}(w^\e)_x+Q^*_{2}(w^\e)_y=-\e\t^\e\nabla(\s_1(\r^\e)\,\nabla\t^\e)
 +\e \int_{\bar{q}}^{q(\r^\e)}\frac{dq}{\r q}\,\,\nabla(\s_2(\r^\e)\nabla\r^\e),
\ee
where
$$
Q_1^*(w)=\r q H^*_{\mu}\cos\t-qH^*_\theta \sin\t, \qquad Q_2^*(w)=\r
q H^*_{\mu}\sin\t + qH^*_\theta \cos\t.
$$
Then \be \label{Q1}
\begin{split}
Q^*_{1}(w^\e)_x+Q^*_{2}(w^\e)_y
 &= \e\div\left(-\s_1(\r^\e)\t^\e\nabla\t^\e
   +\s_2(\r^\e)\big(\int_{\bar{q}}^{q^\e}\frac{dq}{\r
   q}\big)\,\,\nabla\r^\e\right)\\
 &\quad   +\e\s_1(\r^\e)|\nabla\t^\e|^2-\e\s_2(\r^\e)\frac{|\nabla\r^\e|^2}{\r^\e
    q^\e}\frac{dq(\r^\e)}{d\r^\e} \\
 &= \e\div\left(-\s_1(\r^\e)\t^\e\nabla\t^\e
   +\s_2(\r^\e)\nabla\r^\e\int_{\bar{q}}^{q^\e}\frac{dq}{\r
   q}\right)\\
 &\quad   +\e\s_1(\r^\e)|\nabla\t^\e|^2+\e\s_2(\r^\e)\frac{c^2(\r^\e)}{(\r^\e q^\e)^2}
    |\nabla\r^\e|^2.
 \end{split}
 \ee
Integrating \eqref{Q1} over $\O$ and using the divergence theorem,
we obtain with $\bar{q}=u_\infty>0$ that
\begin{equation}  \label{Q6}
\begin{split}
\int_{\d\O}(Q_1^*(w^\e),Q_2^*(w^\e))\cdot{\bf n}\, ds =
 & \e\int_{\d\O}\Big(-\s_1(\r^\e)\t^\e\nabla\t^\e\cdot{\bf n}
   +\s_2(\r^\e)\big(\int_{u_\infty}^{q^\e}\frac{dq}{\r q}\big)\nabla\r^\e\cdot{\bf n}\Big)ds\\
 &+\e\int_\O\left(\s_1(\r^\e)|\nabla\t^\e|^2
 +\s_2(\r^\e)\frac{c^2(\rho^\e)}{(\r^\e q^\e)^2}|\nabla\r^\e|^2\right)dxdy, \\
 :=& I_1+I_2,
\end{split}
\end{equation}
where $Q_1^*,Q_2^*$ depend only on $(\r^\e,\t^\e)$ and are
independent of their derivatives. Thus, from the $L^\infty$ bound,
the left hand side of \eqref{Q6} is uniformly bounded for all
$\e>0$. More specifically, from \eqref{Q} and the formula of $H^*$,
we have the following:

On $\partial\Omega_1$, $(Q_1^*, Q_2^*)\cdot {\bf n}=q^\e\theta^\e$
which is uniformly bounded in $\epsilon>0$;

On the horizontal part of $\partial\Omega_2$, $\theta^\e=0$ and
$|(Q_1^*, Q_2^*)\cdot {\bf n}|=|Q_2^*|=|q^\e\theta^\e|=0$;

On the vertical parts of $\partial\Omega_2$, $\theta^\e=0$ and
$|(Q_1^*, Q_2^*)\cdot {\bf n}|=|Q_1^*|=|\rho(u_\infty) u_\infty
H^*_\mu(u_\infty)|$ is uniformly bounded in $\e>0$.

\medskip
Using the boundary conditions \eqref{bc}, one has
\begin{equation*}
\begin{split}
I_1 =&\e\int_{\d\O_1}\s_2
 \big(\int_{u_\infty}^{q^\e}\frac{dq}{\r q}\big)\nabla\r^\e\cdot{\bf n}\, ds
  +\e\int_{\d\O_2}\s_2
  \big(\int_{u_\infty}^{u_\infty}\frac{dq}{\r q}\big)\,\nabla\r^\e\cdot{\bf n}\,
  ds\\
=&-\e\int_{\d\O_1}|\r^\e(u^\e,v^\e)\cdot{\bf n}|
  \big(\int_{u_\infty}^{q^\e}\frac{dq}{\r q}\big)\, ds,
\end{split}
\end{equation*}
which implies that $I_1$ is uniformly bounded due to the
 $L^\infty$ bound, since the integrand is like $q\ln q$ near $q=0$
and is well behaved. Therefore, \eqref{Q6} implies that $I_2$ is
uniformly bounded.
\end{proof}

In terms of a general generator $H$, equation \eqref{w30} becomes
$$
Q_{1}(w^\e)_x+Q_{2}(w^\e)_y=\e\nabla(\s_1(\r^\e)\nabla\t^\e)(\r^\e
H_{\mu\t}-H_\t)
 +\e\nabla(\s_2(\r^\e)\nabla\r^\e)\big(H_\mu+\frac1\r H_{\t\t}\big),
$$
 or
\be \label{W8}
\begin{split}
Q_{1}(w^\e)_x+Q_{2}(w^\e)_y= \;
 &\e\div\Big(\s_1(\r^\e)\nabla\t^\e (\r^\e H_{\mu\t}-H_\t)
 +\s_2(\r^\e)\nabla\r^\e \big(H_\mu+\frac1\r H_{\t\t}\big)\Big)\\
 &-\e\s_1(\r^\e)\nabla\t^\e\cdot\nabla
 \big(\r^\e H_{\mu\t}-H_\t \big)
 -\e\s_2(\r^\e)\nabla\r^\e\cdot\nabla \big(H_\mu+\frac1{\r^\e} H_{\t\t}\big).
\end{split}
\ee
 The entropy pair $(Q_1,Q_2)$ satisfies \eqref{w31} which can be
 written in terms of $H$,
\be \label{W10}
 \begin{split}
  &\frac{\d Q_1}{\d\r}
   =-\frac{c^2}{\r q}\sin\t\; (\r H_{\mu\t}-H_\t)
    -q\cos\t\; \frac{c^2-q^2}{q^2}\big(H_\mu+\frac1\r H_{\t\t}\big), \\
   &\frac{\d Q_1}{\d\t}
   =q\cos\t\; (\r H_{\mu\t}-H_\t)
   -\r q\sin\t\; \big(H_\mu+\frac1\r H_{\t\t}\big), \\
  &\frac{\d Q_2}{\d\r}
   =\frac{c^2}{\r q}\cos\t\; (\r H_{\mu\t}-H_\t)
    -q\sin\t\; \frac{c^2-q^2}{q^2}\big(H_\mu+\frac1\r H_{\t\t}\big), \\
   &\frac{\d Q_2}{\d\t}
   =q\sin\t\; (\r H_{\mu\t}-H_\t)
   +\r q\cos\t\; \big(H_\mu+\frac1\r H_{\t\t}\big).
   \end{split}
   \ee

\begin{proposition} \label{H-1}
Assume that $q^\e(x,y)\ge \alpha_U$ on $U\subset\Omega$ for some
constant $\alpha_U>0$, in addition $q^\e(x,y)\le q^*<q_{cav}$ensured
by the invariant regions. Then
$$
\partial_xQ_{1\pm}(w^\e)+\partial_y
Q_{2\pm}(w^\e) \quad\mbox{are confined in a compact set in
$H^{-1}(U)$},
$$
for any entropy flux pair $(Q_1,Q_2)$ generated by $H\in C^3$. That
is, hypothesis {\rm (A.2)} of the compensated compactness framework
is satisfied for such entropy pair $(Q_1,Q_2)$ generated by $H\in
C^3$.
\end{proposition}

\begin{proof}
For such an entropy pair $(Q_1,Q_2)$ generated by $H\in C^3$,
equation \eqref{W8} is satisfied. By Proposition \ref{I2}, we have
$$
J_1^\e:=\e\div\Big(\s_1(\r^\e)
 (\r^\e H_{\mu\t}-H_\t)\nabla\t^\e
 +\s_2(\r^\e)\big(H_\mu+\frac1{\r^\e} H_{\t\t}\big)\Big)\nabla\r^\e\to 0
$$
in $H^{-1}(U)$ as $\e\to 0$, and
$$
J_2^\e:=-\e\s_1(\r^\e)\nabla\t^\e\cdot\nabla
 \left(\r^\e H_{\mu\t}-H_\t\right)
 -\e\s_2(\r^\e)\nabla\r^\e\cdot\nabla \big(H_\mu+\frac1{\r^\e} H_{\t\t}\big)
$$
is in $L^1(U)$ uniformly in $\e$. On the other hand, $(Q_1(w^\e),
Q_2(w^\e))$ are uniformly bounded which yield that $J_1^\e+J_2^\e$
is bounded in $W^{-1,\infty}(U)$. Then Murat's lemma \cite{Mu3}
implies that $J_1^\e+J_2^\e$ are confined in a compact subset of
$H^{-1}(U)$.
\end{proof}

As a corollary of Proposition 8.2, we conclude

\begin{proposition}\label{prop8.3}
Let the viscous viscosity fields $w^\e=(u^\e,v^\e)$ have the speed
$q^\e(x,y)$ uniformly bounded in $\e>0$ away from zero when $(x,y)$
is away from the obstacle boundary $\partial\Omega_1$, i.e. there
exists a positive $\e$-independent $\alpha=\alpha(\delta)\to 0$ as
$\delta\to 0$ such that $q^\e(x,y)\ge \alpha(\delta)$ for any
$(x,y)\in \Omega_\delta=\{(x,y)\in\Omega\,:\, dist((x,y),
\partial\Omega_1)\ge \delta>0\}$.
Then
$$
\partial_xQ_{1\pm}(w^\e)+\partial_y
Q_{2\pm}(w^\e) \quad\mbox{are confined in a compact set in
$H^{-1}_{loc}(\Omega)$},
$$
for any entropy pair $(Q_1,Q_2)$ generated by $H\in C^3$.
\end{proposition}

\section{Convergence of the Vanishing Viscosity Solutions}

As noted earlier, the reduction of the support of the Young measure
$\nu$ is accomplished via application of the commutation identity
\eqref{t2}.
Here we follow a technique of Morawetz \cite{Mor95} and use the
entropy generators obtained via classical separation variables first
given by Loewner \cite{Loewner}. Specifically, we look for $H$ as
the following two forms:
$$
H_n(\mu,\t)=F_n(\mu)e^{\pm in\t},
$$
and
$$
H_n(\mu,\t)=K_n(\mu)e^{\pm n\t}.
$$
Hence, from the generalized Tricomi equation \eqref{Tricomi} for
$H$, we see
$$
\ddot{F}_n+n^2\frac{M^2-1}{\r}F_n=0, \qquad
  \ddot{K}_n-n^2\frac{M^2-1}{\r}K_n=0,
$$
where $\dot{}=\frac{d}{d\mu}$. Substitution into the representation
for $Q_1, Q_2$ as given in Proposition \ref{Q8} generates an
infinite sequence of entropy pairs $(Q^{(n)}_{1\pm},
Q^{(n)}_{2\pm})$ associated with $F_n$; and a similar construction
can be done for $K_n$.

First we apply the commutation identity to $(Q^{(n)}_{1\pm},
Q^{(n)}_{2\pm})$ associated with $F_n$. Set
$$
I=(Q_{1+}^{(n)}-{Q_{1+}^{(n)}}')(Q_{2-}^{(n)}-{Q_{2-}^{(n)}}')
-(Q_{2+}-{Q_{2+}^{(n)}}')(Q_{1-}-{Q_{1-}^{(n)}}').
$$
Then
\begin{equation*}
\begin{split}
I=&\left((\dot{F}_{n}\r q\cos\t + in{F}_{n}q\sin\t)e^{in\t}
-(\dot{F}'_{n}\r' q'\cos\t' +in{F}'_{n}q'\sin\t')e^{in\t'}
 \right) \\
&\,\,\,\times\left((\dot{F}_{n}\r q\sin\t +
in{F}_{n}q\cos\t)e^{-in\t} -(\dot{F}'_{n}\r' q'\sin\t'
+in{F}'_{n}q'\cos\t')e^{-in\t'}
\right) \\
&-\left((\dot{F}_{n}\r q\sin\t - in{F}_{n}q\cos\t)e^{in\t}
-(\dot{F}'_{n}\r' q'\sin\t' - in{F}'_{n}q'\cos\t')e^{in\t'}
\right) \\
&\quad\times\left((\dot{F}_{n}\r q\cos\t -
in{F}_{n}q\sin\t)e^{-in\t} -(\dot{F}'_{n}\r' q'\cos\t' -
in{F}'_{n}q'\sin\t')e^{-in\t'}
  \right).
\end{split}
\end{equation*}
With a tedious calculation, we obtain
\begin{equation*}
\begin{split}
&<\nu\otimes\nu', \quad I> \\
&=  \Big< \nu\otimes\nu', \quad
  2in\r q^2\left({F}_{n}\dot{F}_{n}+{F}_{n}\dot{F}_{n}\right)\\
 &\qquad\qquad\qquad +2e^{in(\t-\t')}\left(-\dot{F}_{n}\r q\cos\t\dot{F}'_{n}\r'q'\sin\t'
 +n^2{F}_{n}q\sin\t{F}'_{n}q'\cos\t'\right.\\
 &\qquad\qquad\qquad\qquad\qquad\qquad\left.+\dot{F}_{n}\r q\sin\t\dot{F}'_{n}\r'q'\cos\t'
 -n^2{F}_{n}q\cos\t{F}'_{n}q'\sin\t'\right)\\
 &\qquad\qquad\qquad +2e^{in(\t-\t')}in\left(-\dot{F}_{n}\r q\cos\t\dot{F}'_{n}q'\cos\t'
 -{F}_{n}q\sin\t{F}'_{n}\r'q'\sin\t'\right.\\
 &\qquad\qquad\qquad\qquad\qquad\qquad\,\,\left.-\dot{F}_{n}\r q\sin\t\dot{F}'_{n}q'\sin\t'
 -{F}_{n}q\cos\t{F}'_{n}\r'q'\sin\t'\right)\Big>,
\end{split}
\end{equation*}
and then
\begin{equation*}
\begin{split}
 &\Big<\nu\otimes\nu', \quad -\frac{i}{4}I\Big> \\
 &= \Big<\nu\otimes\nu', \quad
   n\r q^2F_n\dot{F}_n
  +\frac12F_nF'_nqq'\sin(n(\t-\t'))\sin(\t-\t')(\r\r'+n^2)\\
 & \qquad\qquad\qquad -n\dot{F}_nF'_n\r qq'\cos(n(\t-\t'))\cos(\t-\t')\Big>\\
 &= \Big<\nu\otimes\nu', \quad \frac{n}2(qF_n-q'F'_n)(\r q\dot{F}_n-\r'q'\dot{F}'_n)\\
 & \qquad\qquad\qquad +\frac12\sin^2\left(\frac{n+1}2(\t-\t')\right)
   \left(2n\dot{F}_n{F}'_n\r qq'+F_nF'_nqq'(\r\r'+n^2)\right)\\
 & \qquad\qquad\qquad +\frac12\sin^2\left(\frac{n-1}2(\t-\t')\right)
   \left(2n\dot{F}_n{F}'_n\r qq'-F_nF'_nqq'(\r\r'+n^2)\right)\Big>.
\end{split}
\end{equation*}
Thus, from \eqref{t2},
\begin{equation*}
\begin{split}
<\nu\otimes\nu', \quad
 & \frac{n}2(qF_n-q'F'_n)(\r q\dot{F}_n-\r'q'\dot{F}'_n)\\
 & +\frac12\sin^2\left(\frac{n+1}2(\t-\t')\right)
   \left(2n\dot{F}_n{F}'_n\r qq'+F_nF'_nqq'(\r\r'+n^2)\right)\\
 & +\frac12\sin^2\left(\frac{n-1}2(\t-\t')\right)
   \left(2n\dot{F}_n{F}'_n\r qq'-F_nF'_nqq'(\r\r'+n^2)\right)>=0.
\end{split}
\end{equation*}
Now interchange $n$ and $-n$ to get the $n^2$-terms equal to zero.
Hence
\begin{equation} \label{young}
\begin{split}
<\nu\otimes\nu', \quad
 & (qF_n-q'F'_n)(\r q\dot{F}_n-\r'q'\dot{F}'_n)\\
 & +2\left(\sin^2\left(\frac{n+1}2(\t-\t')\right)
   +\sin^2\left(\frac{n-1}2(\t-\t')\right)\right)\dot{F}_n{F}'_n\r qq'>=0.
\end{split}
\end{equation}

Then we follow verbatim from the argument of Morawetz \cite{Mor95}
to lead the following convergence theorem.

\begin{theorem}
Let $v_\infty=0, \, |u_\infty|< q_{cav}, \, 1\le\g<3$. Assume that
there exists a positive $\e$-independent $\alpha(\delta)\to 0$ as
$\delta\to 0$ such that $q^\e(x,y)\ge \alpha(\delta)$ for any
$(x,y)\in \Omega_\delta=\{(x,y)\in\Omega\,:\, dist((x,y),
\partial\Omega_1)\ge \delta>0\}$. Then the following hold:
\begin{enumerate}
\item[{\rm (i)}] The support of the Young measure $\nu_{x,y}$ strictly excludes
the stagnation point $q=0$ and reduces to a point for a.e.
$(x,y)\in\O_\delta$, hence the Young measure is a Dirac mass;

\item[{\rm (ii)}] The sequence $(u^\e,v^\e)$ has a  subsequence
converging strongly in $L^2_{loc}(\O)^2$  to an entropy solution in
the sense of Definition {\rm 7.1};

\item[\rm (iii)] The boundary condition $(u,v)\cdot{\bf n}\ge 0$ on
$\d\O_1$ is satisfied in the sense of normal trace in \cite{CFrid}.
%
\end{enumerate}
\end{theorem}

\begin{proof}
First, we choose $F_n\sim q^{-n}$, with $n$ sufficiently large and
$q<q_{cr}$, to show that the support of $\nu$ cannot lie in both
$q<q_{cr}$ and $q>q_{cr}$ as in Morawetz \cite{Mor85}. If the
support of $\nu$ lies inside $q\le q_{cr}$, equation \eqref{young}
or the argument of Chen-Dafermos-Slemrod-Wang \cite{CDSW} applies to
show that $\nu$ is a Dirac mass. If the support of $\nu$ is in $q\ge
q_{cr}$, Morawetz's reconstruction of DiPerna's argument in
\cite{Mor85} by using $K_n$ associated entropies $(Q_{1\pm},
Q_{2\pm})$ works to show again that $\nu$ is a Dirac mass. The
boundary condition on $\d\O_1$ follows the weak form of the second
equation in \eqref{vcl} with
$R_2=\e\nabla(\sigma_2(\rho)\nabla\rho)$. The entropy inequality
follows from \eqref{Q1}.
\end{proof}

\begin{remark}\label{rm9.1}
From the boundary condition \eqref{bc} for the viscous problem, we
find that, if one can show that $\e
\sigma_2(\rho^\e)\nabla\rho^\e\cdot {\bf n}\to 0$ in the weak sense
on $\partial\Omega_1$ as $\e\to 0$, then we can conclude the desired
boundary condition $(u,v)\cdot{\bf n}= 0$ on $\d\O_1$ satisfied in
the weak sense. This is the case when the limit solution has certain
regularity along the boundary, i.e. any shock strength is zero at
its intersection point with the boundary, as conjectured for the
shock formed in the supersonic bubble near the obstacle (cf.
Morawetz \cite{Mor82,Mor04}).
\end{remark}

\begin{remark}\label{rm9.2}
For the first boundary value problem in Figure 6(a), if
$\partial\Omega_1$ is sufficiently smooth, one expects that the
transonic flow is strictly away from stagnation, which would expect
that the viscous solutions stay uniformly away from cavitation in
the whole domain $\Omega$. However, for the second boundary value
problem in Figure 6(b), one expects that the transonic flow
generally has a stagnation point on the boundary near the head of
the obstacle since the circulation about the boundary
$\partial\Omega_2$ is zero. The assumption that there exists an
$\e$-independent positive $\alpha(\delta)\to 0$ as $\delta\to 0$
such that $q^\e(x,y)\ge \alpha(\delta)$, for any $(x,y)\in
\Omega_\delta$, is physically motivated and designed in order to fit
both cases.
\end{remark}

\section{Resolution of the Viscous Problem}

As we have seen in the previous sections, it is crucial that our
viscous problem possesses smooth solutions $w^\e=(u^\e,v^\e)$. We
address the issue in this section.
Consider the viscous problem:
\begin{equation} \label{vc2}
\begin{cases}
v_x-u_y=\e\nabla(\s_1(\r)\nabla\t), \\
(\r u)_x+(\r v)_y=\e\nabla(\s_2(\r)\nabla\r),
\end{cases}
\end{equation}
with
\begin{equation} \label{sigma}
 \s_1=1, \quad
 \s_2(\r)=\frac{q^2-c^2}{q^2} \qquad\text{for}\quad
 q\ge\sqrt{2}\,q_{cr},
\end{equation}
and $\s_1, \s_2$ any smooth bounded, positive continuation for
$q\le\sqrt{2}q_{cr}$, and with the boundary conditions on the
obstacle $\d\O$:
\begin{equation} \label{bc1}
\begin{cases}
\nabla\t\cdot{\bf n}=0 \qquad\text{on}\quad \d\O_1,\\
\e\s_2\nabla\r\cdot{\bf n}-|\r\;(u,v)\cdot{\bf n}|=0
  \qquad\text{on}\quad \d\O_1, \\
(u,v)-(u_\infty,v_\infty)=0  \qquad\text{on}\quad \d\O_2,
 \quad\text{with}\quad q_\infty< q_{cr},
\end{cases}
\end{equation}
where $\r$ is a function of $q$ given by Bernoulli's law
\eqref{ber1} for $\g>1$ and \eqref{ber2} for $\g=1$, and ${\bf n}$
is the unit normal pointing into the flow region on $\d\O$.

Introduce a new variable
$$
\s(\r):=\int_1^\r\s_2(\xi)d\xi.
$$
Since $\s'(\r)=\s_2(\r)>0$, we can always invert to write
 $\r=\r(\s)$, as well as $q=q(\s)$.
 The advantage of this change of dependent variable is that we now
 have the viscous terms to be $\e\Delta\t$ and $\e\Delta\s$.

 We use the polar velocity notation
 $$u=q(\s)\cos\t, \quad v=q(\s)\sin\t,$$
 and write the viscous system \eqref{vc2} as
\begin{equation} \label{vc3}
\begin{cases}
(q(\s)\sin\t)_x-(q(\s)\cos\t)_y=\e\Delta\t, \\
(\r(\s)q(\s)\cos\t)_x+(\r(\s)q(\s)\sin\t)_y=\e\Delta\s,
\end{cases}
\end{equation}
which yields
\begin{equation} \label{vc4}
\begin{cases}
q'(\s)\sin\t \, \s_x +q(\s)\cos\t \, \t_x
-q'(\s)\cos\t \, \s_y +q(\s)\sin\t \, \t_y=\e\Delta\t, \\
(\r(\s)q(\s))'\cos\t \, \s_x -\r(\s)q(\s)\sin\t \, \t_x
+(\r(\s)q(\s))'\sin\t \, \s_y +\r(\s)q(\s)\cos\t \,
\t_y=\e\Delta\s,
\end{cases}
\end{equation}
with the boundary conditions on $\d\O_1$:
\begin{equation} \label{bc2}
\begin{cases}
\nabla\t\cdot {\bf n}=0, \\
\e\nabla\s\cdot {\bf n}-|\r(\s)q(\s)(\cos\t, \sin\t)\cdot {\bf
n}|=0,
\end{cases}
\end{equation}
and the boundary conditions on $\d\O_2$:
\begin{equation} \label{bc21}
\s = \s_\infty, \quad \t =\t_\infty,
\end{equation}
where $\s_\infty=\s(\r_\infty)$ and $\t_\infty$ are constants. It is
convenient to have the homogeneous boundary conditions on $\d\O_2$:
$$\bar\s:=\s-\s_\infty, \quad \bar\t:=\t-\t_\infty; \qquad
\bar q(\bar\s):=q(\bar\s+\s_\infty), \quad
\bar\r(\bar\s):=\r(\bar\s+\s_\infty).$$
 Hence, we have the following boundary value problem:
\begin{equation} \label{vc5}
\begin{cases}
\e\Delta\bar\t \, =
 & \bar q'(\bar\s)\sin(\bar\t+\t_\infty) \,\bar\s_x
      +\bar q(\bar\s)\cos(\bar\t+\t_\infty) \, \bar\t_x \\
 &-\bar q'(\bar\s)\cos(\bar\t+\t_\infty) \, \bar\s_y
      +\bar q(\bar\s)\sin(\bar\t+\t_\infty) \, \bar\t_y, \\
\e\Delta\bar\s \, =
 &(\bar\r(\bar\s)\bar q(\bar\s))'\cos(\bar\t+\t_\infty) \, \bar\s_x
     -\bar\r(\bar\s)\bar q(\bar\s)\sin(\bar\t+\t_\infty) \, \bar\t_x \\
 &+(\bar\r(\bar\s)\bar q(\bar\s))'\sin(\bar\t+\t_\infty) \, \bar\s_y
     +\bar\r(\bar\s)\bar q(\bar\s)\cos(\bar\t+\t_\infty) \, \bar\t_y,
\end{cases}
\end{equation}
with the boundary conditions on $\d\O_1$:
\begin{equation} \label{bc3}
\begin{cases}
\e\nabla\bar\t\cdot {\bf n}=0, \\
\e\nabla\bar\s\cdot {\bf n}
  -|\bar\r(\bar\s)\bar q(\bar\s)(\cos(\bar\t+\t_\infty),
       \sin(\bar\t+\t_\infty))\cdot {\bf  n}|=0,
\end{cases}
\end{equation}
and $(\bar\s, \bar\t)=(0,0)$ on $\d\O_2$.

Consider the classical equation
$$(\bar\t,\bar\s)=\l\Gamma(\bar\t,\bar\s), \quad 0\le\l\le 1,$$
which corresponds to the boundary value problem:
\begin{equation} \label{vc7}
\begin{split}
&\e\Delta\bar\t =
  \l\bar q'(\bar\s)\sin(\bar\t+\t_\infty) \, \bar\s_x
      +\l\bar q(\bar\s)\cos(\bar\t+\t_\infty) \, \bar\t_x \\
 &\qquad\qquad -\l\bar q'(\bar\s)\cos(\bar\t+\t_\infty) \, \bar\s_y
      +\l\bar q(\bar\s)\sin(\bar\t+\t_\infty) \, \bar\t_y
      \qquad\text{in}\,\,\, \O, \\
&\e\Delta\bar\s =
 \l(\bar\r(\bar\s)\bar q(\bar\s))'\cos(\bar\t+\t_\infty)
    \, \bar\s_x
     -\l\bar\r(\bar\s)\bar q(\bar\s)\sin(\bar\t+\t_\infty)
      \, \bar\t_x \\
 &\qquad\qquad +\l(\bar\r(\bar\s)\bar q(\bar\s))'
 \sin(\bar\t+\t_\infty)    \, \bar\s_y
     +\l\bar\r(\bar\s)\bar q(\bar\s)\cos(\bar\t+\t_\infty)
     \, \bar\t_y \qquad\text{in}\,\,\, \O,\\
&\e\nabla\bar\t\cdot {\bf n} = 0
   \qquad\text{on}\,\,\, \d\O_1,\\
&\e\nabla\bar\s\cdot{\bf n} =
  \l|\bar\r(\bar\s)\bar q(\bar\s)(\cos(\bar\t+\t_\infty),
      \sin(\bar\t+\t_\infty))\cdot {\bf n}| \qquad\text{on}\,\,\, \d\O_1, \\
&(\bar\s, \bar\t)=(0,0) \qquad\text{on}\,\,\, \d\O_2,
\end{split}
\end{equation}
which is the original system with $\e$ replaced by $\e\l^{-1}$.

Now we recall some results from the linear elliptic theory. Consider
the mixed Dirichlet-Neumann problem
\begin{equation} \label{e1}
\begin{cases}
\Delta w=f_1 &\text{in}\,\,\, \O, \\
\frac{\d w}{\d{\bf n}}=f_2 &\text{on}\,\,\, \d\O_1,\\
w=f_3 &\text{on}\,\,\, \d\O_2,
\end{cases}
\end{equation}
where $\O$ is a bounded domain in $\R^2$ with the boundary
$\d\O=\d\O_1\cup\d\O_2$. We now introduce weighted H\"{o}lder norms.
For an arbitrary open set $S$ and an arbitrary $a>0$, we have the
usual H\"{o}lder norms $\|\cdot\|_{a,S}$. If $\delta>0$ and $a+b\ge
0$, we define
 $$\O_\delta=
  \left\{x\in\O:\; {\rm dist}(x,\d\O_2)>\delta\right\},$$
and
$$\|w\|_a^{(b)}=\sup_{\delta>0}\delta^{a+b}\|w\|_{a,\O_\delta},$$
and then denote $C^a_{(b)}$ the set of all functions with
$\|w\|_a^{(b)}<\infty$, and $C^*(\O)=C^2(\O\cup\d\O_1)\cap
C(\bar\O)$.

\begin{lemma}[Lieberman\cite{Lieberman}, Theorem 2(b)]
  Suppose $\d\O_1\in C^{2+\alpha}$, $\d\O_2$ is nonempty, and $\O$
  satisfies the global $\Sigma$-wedge condition and the uniform
  interior and exterior cone condition on $\d\O$. Then there is a
  unique solution to the linear elliptic problem \eqref{e1} and
$$
\|w\|_{2+\alpha}^{(\beta)}\le
  C\left(\|f_1\|_\a^{(2+\beta)}+\|f_2\|_{1+\a}^{(1+\beta)}
  +\|f_3\|_\beta\right)
$$
for  some $\beta>0$ and constant $C>0$.
\end{lemma}

\begin{remark}
Note that Lieberman's theorem is given for the boundary condition on
$\d\O_1$: $Mw:=\beta^iD_i w+d\, w=f_2$, where $d<0$. However, as
noted by Lieberman on page 429 in \cite{Lieberman}, the Fredholm
alternative holds, and if $\d\O_2$ is nonempty, the restriction
$d<0$ may be relaxed to $d\le 0$. The domains shown in Figure
\ref{domain} satisfy the above conditions.
\end{remark}

  We define $\Gamma: (\bar\Theta,\bar\Sigma)\to(\bar\t,\bar\s)$ to
  be the solution map associated with solving the decoupled mixed
  Dirichlet-Neumann boundary value problem obtained from
  substitution of $(\bar\Theta,\bar\Sigma)$ into the right hand
  side of \eqref{vc7}. For $(\bar\Theta,\bar\Sigma) \in
  (C^{1+\a}_{(2+\beta)})^2$, the associated functions $f_1$ and $f_2$
  give finite values on the right-hand side of the
 Schauder estimates of Lieberman's theorem and $f_3=0$ in our
  problem. Therefore, $(\bar\t,\bar\s)$ stays in a bounded subset of
  $(C^{2+\a}_{(\beta)})^2$ which is a compact subset of
  $(C^{1+\a}_{(2+\beta)})^2$, and $\Gamma$ is a compact map of
  the Banach space $B=(C^{1+\a}_{(2+\beta)})^2$ into itself.

We recall one popular version of the Leray-Schauder fixed point
theorem.
\begin{lemma}[Leray-Schauder Fixed Point Theorem]
Let $\Gamma$ be a compact mapping of a Banach space $B$ into
itself, and suppose that there exists a constant $M$ such that
$\|W\|_{B}\le M$ for all $W\in B$ and all $\l\in[0,1]$ satisfying
$W=\l \, \Gamma W$. Then $\Gamma$ has a fixed point.
\end{lemma}

To apply the above theorem to our example, we first note that the
only solution when $\l=0$ is $(\bar\t,\bar\s)=(0,0)$, hence it is
certainly bounded in $B$. For $\l\neq 0$, we study the solutions to
our original viscous system with $\bar\e=\e\l^{-1}\ge\e$. The
invariant region argument we have given provides the uniform
$L^\infty(\bar\O)$ bound independent of $\bar\e$. The gradient
estimates from the entropy equality yield that
$\sqrt{\bar\e}\nabla(\bar\t,\bar\s)$ are in a bounded set of
$(L^2(\O))^2$, and hence $\sqrt{\e}\nabla(\bar\t,\bar\s)$ are in the
same bounded set of $(L^2(\O))^2$ for all $0<\l\le 1$.

Now we examine the equations themselves:
$$\e\Delta\bar\t=\l f_1, \qquad \e\Delta\bar\s=\l f_2,$$
where $f_1, f_2$ are in a bounded set of $L^2(\O)$. Hence, for
$0<\l\le 1$, $\e\Delta(\bar\t, \bar\s)$ are in a bounded set of
$L^2(\O)$ and, by the Calderon-Zygmund inequality (Theorem 9.9 in
Gilbarg-Trudinger \cite{GT}), $(\bar\t,\bar\s)\in W^{2,2}(\O)$ .

Next we differentiate the both sides of our equations with respect
to $x, y$, and we get the equations of the form
$$
\e\Delta\bar\t_x=\l f_{1x}, \qquad \e\Delta\bar\t_y=\l f_{1y}.
$$
Since the right-hand sides are in a bounded set of $L^2(\O)$, so are
the left-hand sides. We continue in this fashion to see that there
is a constant $M_1$ (independent of $m$ and $\l$) such that
$$
\|(\bar\t,\bar\s)\|_\infty+
\|\nabla(\bar\t,\bar\s)\|_{(H^m(\O))^2}\le M_1,
$$
for all $1\le m<\infty$.

Then Morrey's embedding theorem applied to our case yields
$\nabla(\bar\t,\bar\s)$ bounded by constant $M_1=M$ in
$(C^s(\bar\O))^2$ for any $s\ge 1$ and $(\bar\t, \bar\s)$ bounded in
$B$ for any $\tau>0$ by the constant $M_1=M$, where the constant is
independent of $\l$. Thus, all possible solutions of
$(\bar\t,\bar\s)=\l \, \Gamma(\bar\t,\bar\s)$, $0\le\l\le 1$, are
bounded by constant $M_1=M$ with $M$ independent of $\l$. Therefore,
the hypotheses of the Leray-Schauder fixed point theorem are
satisfied and the viscous problem has a solution, i.e. we have the
following theorem.

\begin{theorem}
The viscous problem has a solution $(u^\e, v^e)$ for every $\e>0$
such that $q^\e\le q^*<q_{cav}$ for $\gamma\in [1, 3)$, where $q^*$
is a constant independent of $\e>0$.
\end{theorem}

\bigskip
\bigskip
\noindent{\bf Acknowledgments.} Gui-Qiang Chen's research was
supported in part by the National Science Foundation under Grants
DMS-0505473, DMS-0244473, and an Alexandre von Humboldt Foundation
Fellowship.
Marshall Slemrod's research was supported in part by the National
Science Foundation grant DMS-0243722.
Dehua Wang's research was supported in part by the National
Science Foundation grant DMS-0244487 and the Office of Naval
Research grant N00014-01-1-0446.
The authors would like to thank the other members of the NSF-FRG on
multidimensional conservation laws: S. Canic, C.~D. Dafermos, J.
Hunter, T.-P. Liu, C.-W. Shu, and Y. Zheng for stimulating
discussions.

\end{document}